\magnification=1000
\hsize=11.7cm
\vsize=18.9cm
\lineskip2pt \lineskiplimit2pt
\nopagenumbers

\hoffset=-1truein
\voffset=-1truein

\advance\voffset by 4truecm
\advance\hoffset by 4.5truecm

\newif\ifentete

\headline{\ifentete\ifodd	\count0 
      \rlap{\head}\hfill\tenrm\llap{\the\count0}\relax
    \else
        \tenrm\rlap{\the\count0}\hfill\llap{\head} \relax
    \fi\else
\global\entetetrue\fi}

\def\entete#1{\entetefalse\gdef\head{#1}} 
\entete{}

\input amssym.def
\input amssym.tex

\def\-{\hbox{-}}
\def\.{{\cdot}}
\def\O{{\cal O}}
\def\K{{\cal K}}
\def\F{{\cal F}}

\def\P{{\cal P}}

\def\T{{\cal T}}

\def\X{{\cal X}}

\def\ch{\frak c\frak h}

\def\Gr{\frak G\frak r}

\def\id{\frak i\frak d}
\def\int{\frak i\frak n\frak t}

\def\qq{\quad{\rm and}\quad}

\def\too{\longrightarrow}
\def\aut{\frak a\frak u\frak t}

\def\Loc{\frak L\frak o\frak c}
\def\loc{\frak l\frak o\frak c}

 3
 2
\font\large=cmr10  scaled \magstep 2
 2
 2
 1
 2

\font\cds=cmr7
\font\cdt=cmti7

\count0=1

\centerline{\large  Beyond a question of Markus Linckelmann}
\medskip
\centerline{\bf Lluis Puig }
\medskip
\noindent 
\centerline{\cds CNRS, Institut de Math\'ematiques de Jussieu, lluis.puig@imj-prg.fr}
\par
\noindent
\centerline{\cds 6 Av Bizet, 94340 Joinville-le-Pont, France}

\medskip
\noindent
{\bf Abstract:} {\cds  In the 2002 Durham Symposium, Markus Linckelmann conjectured the existence of a {\cdt regular central ${\scriptstyle k^*}\-$extension\/}  of the full subcategory over the 
{\cdt selfcentralizing Brauer pairs\/} of the {\cdt Frobenius ${\scriptstyle P}\-$category\/} 
${\scriptstyle \F_{(b,G)}}$ associated with a block ${\scriptstyle b}$  of defect group ${\scriptstyle P}$ of a finite group ${\scriptstyle G}\,,$ which would include,
as ${\scriptstyle k^*}\-$automorphism groups of the objects, the ${\scriptstyle k^*}\-$groups associated with the {\cdt automizers\/} of the corresponding {\cdt selfcentralizing Brauer pairs\/}, introduced in [4,~6.6].  As a matter of fact, in this question the {\cdt selfcentralizing Brauer pairs\/} can be replaced by the {\cdt nilcentralized Brauer pairs\/}, still getting a positive answer. But the condition on the  
${\scriptstyle k^*}\-$automorphism groups of the objects is {\cdt not\/} precise enough to guarantee the {\cdt uniqueness\/} of a solution, as showed in [3,~Theorem~1.3]. This {\cdt uniqueness\/} depends on the {\cdt folder structure\/} [6,~Section~2] associated with ${\scriptstyle \F_{(b,G)}}$ in 
[5,~Theorem~11.32], and here we prove the {\cdt existence\/} and the {\cdt uniqueness\/} for any
{\cdt folded Frobenius ${\scriptstyle P}\-$category\/}\/}.

\bigskip
\noindent
{\bf £1. Introduction }
\medskip

£1.1. Let $p$ be a prime number and $\O$ a complete discrete valuation ring with a 
{\it field of quotients\/} $\K$ of characteristic zero and a {\it residue field\/} $k$ of characteristic $p\,;$ we assume that $k$ is algebraically closed. Let $G$ be a finite group, $b$ a {\it block\/} of $G$ --- namely a primitive idempotent in the center $Z(\O G)$ of the group $\O\-$algebra $\O G$ --- and $(P,e)$ a maximal {\it Brauer $(b,G)\-$pair\/} [5,~1.16]; recall that the {\it Frobenius $P\-$category\/} $\F_{\!(b,G)}$ associated with~$b$ is the subcategory of the category of finite groups where the objects are all the subgroups of $P$ and, for any pair of subgroups $Q$ and $R$ of~$P\,,$ the morphisms~$\varphi$ from $R$ to $Q$ are the group homomorphisms $\varphi\,\colon R\to Q$ induced by the conjugation of some element $x\in G$ fulfilling
$$(R,g)\i (Q,f)^x
\eqno £1.1.1\phantom{.}$$
where $(Q,f)$ and $(R,g)$ are the corresponding Brauer $(b,G)\-$pairs contained in 
$(P,e)$ [5, Ch.~3]. 

\medskip
£1.2. Moreover, we say that a Brauer $(b,G)\-$pair $(Q,f)$ is {\it nilcentralized\/} if~$f$ is a {\it nilpotent block\/} of $C_G(Q)$ [5,~7.4], and that $(Q,f)$ is {\it selfcentralizing\/} if the image $\bar f$ of $f$ is a block of {\it defect zero\/} of  $\bar C_G (Q) = C_G(Q)/Z(Q)$~[5,~7.4]; thus, a selfcentralizing Brauer $(b,G)\-$pair is still nilcentralized.  We respectively denote by $\F_{\!(b,G)}^{^{\rm nc}}$ or $\F_{\!(b,G)}^{^{\rm sc}}$ the {\it full\/} subcategories of $\F_{\!(b,G)}$ over the set of subgroups $Q$ of $P$ such that the Brauer $(b,G)\-$pair $(Q,f)$ contained in~$(P,e)$ is respectively nilcentralized or selfcentraling.

\medskip
£1.3. Recall that a {\it $k^*\-$group\/} $\hat G$ is a group endowed with an injective group homomorphism $\theta\,\colon k^*\to Z(\hat G)$~[4,~\S5], that $G = \hat G/\theta (k^*)$ is called the {\it $k^*\-$quotient\/} of~$\hat G$  and that a {\it $k^*\-$group homomorphism\/} is a group
homomorphism which preserves the ``multiplication'' by $k^*\,;$ let us denote by $k^*\-\Gr$ the category of $k^*\-$groups with finite $k^*\-$quotient. In the case of the {\it Frobenius $P\-$category\/} 
$\F_{\!(b,G)}$ above, for any {\it nilcentralized\/} Brauer $(b,G)\-$pair $(Q,f)$ contained in~$(P,e)$ it is well-kown that the action of $N_G (Q,f)$ on the
simple algebra $\O C_G(Q) f/J\big(\O C_G(Q) f\big)$ suplies a $k^*\-$group
$\hat N_G(Q,f)/C_G (Q)$ of $k^*\-$quotient $\F_{\!(b,G)}(Q)\cong N_G(Q,f)/C_G (Q)$ [5,~7.4].

\medskip
£1.4. On the other hand, for any category $\frak C$ and any
Abelian group~$Z$ let us call {\it regular central $Z\-$extension\/} of $\frak C$ any category $\hat\frak C$ over the same objects endowed with a {\it full\/} functor $\frak c\,\colon \hat\frak C\to \frak C\,,$ which is the identity over the objects, and, for any pair of $\frak C\-$objects~$A$ and $B\,,$ with a 
{\it regular\/} action of $Z$  over the {\it fibers\/} of the map
$$\hat\frak C (B,A)\too \frak C (B,A)
\eqno £1.4.1\phantom{.}$$
induced by $\frak c$ --- where $\frak C (B,A)$ and $\hat\frak C (B,A)$ denote the corresponding  sets of 
$\frak C\-$ and $\hat\frak C\-$morphisms from $A$ to $B$ --- in  such a way that these $Z\-$actions are compatible with the composition of $\hat\frak C\-$morphisms. Note that, if $\,\frak C'$ is a second category and $\frak e\, \,\colon  \frak C\to \frak C'$ an equivalence of categories, we easily can obtain 
a {\it regular central $Z\-$extension\/} $\hat\frak C'$ of~$\frak C'$ and a 
{\it $Z\-$compatible equivalence of categories\/}  $\hat\frak e\,\colon \hat\frak C\to  \hat\frak C'\,.$ 
In short, we call {\it $k^*\-$category\/} any {\it regular central $k^*\-$extension\/} of a category.

\medskip
£1.5. In the 2002 Durham Symposium, Markus Linckelmann conjectured
 the existence of a {\it regular central $k^*\-$extension\/} 
 $\hat\F_{\!(b,G)}^{^{\rm sc}}$ of $\F_{\!(b,G)}^{^{\rm sc}}$ admitting
 a {\it $k^*\-$group isomorphism\/}
$$\hat\F_{\!(b,G)}^{^{\rm sc}} (Q)  \cong
\hat N_G(Q,f)/C_G (Q)
\eqno £1.5.1$$
for any {\it selfcentralizing\/}  Brauer $(b,G)\-$pair $(Q,f)$ contained in~$(P,e)\,.$
Here we show the existence of a {\it regular central $k^*\-$extension\/} 
 $\hat\F_{\!(b,G)}^{^{\rm nc}}$ of $\F_{\!(b,G)}^{^{\rm nc}}$ admitting
 a {\it $k^*\-$group isomorphism\/}
$$\hat\F_{\!(b,G)}^{^{\rm nc}} (Q)  \cong
\hat N_G(Q,f)/C_G (Q)
\eqno £1.5.2$$
for any {\it nilcentralized\/}  Brauer $(b,G)\-$pair $(Q,f)$ contained in~$(P,e)\,,$
proving Linckelmann's conjecture.

\medskip
£1.6. In both cases, these {\it $k^*\-$group isomorphisms\/} are not precise enough
to guarantee the uniqueness either of $\hat\F_{\!(b,G)}^{^{\rm nc}}\,,$ or of 
$\hat\F_{\!(b,G)}^{^{\rm sc}}$ as showed in [3,~Theorem~1.3]. More explicitly, if
$(Q,f)$ and $(R,g)$ are {\it nilcentralized\/}  Brauer $(b,G)\-$pairs  contained in~$(P,e)$
such that $(R,g)$ is contained and normal in~$(Q,f)$ then, denoting by  $\hat N_G (Q,f)_R$
the stabilizer of $R$ in $\hat N_G (Q,f)\,,$ Proposition~11.23 in [5] suplies a particular
$k^*\-$group homomorphism
$$\hat N_G(Q,f)_R/C_G (Q)\too \hat N_G(R,g)/C_G (R)
\eqno £1.6.1.$$
But,  a {\it regular central $k^*\-$extension\/}  $\hat\F_{\!(b,G)}^{^{\rm nc}}$ of 
$\F_{\!(b,G)}^{^{\rm nc}}$ also suplies a $k^*\-$group homomorphism
$$\hat\F_{\!(b,G)}^{^{\rm nc}} (Q)_R\too \hat\F_{\!(b,G)}^{^{\rm nc}} (R)
\eqno £1.6.2,$$
where $\hat\F_{\!(b,G)}^{^{\rm nc}} (Q)_R$ denotes the stabilizer of $R$ in 
$\hat\F_{\!(b,G)}^{^{\rm nc}} (Q)\,,$ sending any $\hat\sigma$ 
in~$\hat\F_{\!(b,G)}^{^{\rm nc}} (Q)_R$ on the unique element 
$\hat \tau \in \hat\F_{\!(b,G)}^{^{\rm nc}} (R)$ fulfilling $\hat\iota_R^Q\circ \hat\tau =
\hat\sigma\circ \hat\iota_R^Q\,,$ where $\hat\iota_R^Q$ is a lifting to 
$\hat\F_{\!(b,G)}^{^{\rm nc}} (Q,R)$ of the inclusion map $R\i Q\,.$
The~{\it uniqueness\/} of a suitable {\it regular central $k^*\-$extension\/}  
$\hat\F_{\!(b,G)}^{^{\rm nc}}$ depends on the com-patibility 
of all the $k^*\-$group homomorphisms~£1.6.1 and £1.6.2 with the cor-responding $k^*\-$group isomorphisms~£1.5.2 or, more generally, it depends on the {\it folded structure\/} of 
$\F_{\!(b,G)}^{^{\rm nc}}$ determined by [5,~Theorem~11.32].

 \bigskip
\noindent
{\bf £2.  Folded Frobenius $P\-$categories\/}

\medskip 
£2.1. Denoting by  $P$ a finite $p\-$group, by  $\frak i\Gr$ the category formed by the finite groups and by the injective group  homomorphisms, and  by $\F_{P}$ the subcategory of  $\frak i\Gr$ where the objects are all the  subgroups of $P$ and the morphisms are the group homomorphisms induced by the conjugation by elements of $P\,,$ recall that  a {\it Frobenius  $P\-$category\/} $\F$ is a subcategory 
of $\frak i\Gr$ containing $\F_{P}$ where the objects are all the  subgroups of $P$
and the morphisms fulfill the following three conditions [5, 2.8 and Proposition~2.11]
\smallskip
\noindent
£2.1.1\quad {\it If $Q\,,$ $R$ and $T$ are subgroups of $P\,,$ for any
$\varphi\in \F (Q,R)$ and any group homomorphism $\psi\,\colon T\to R$
the composition $\varphi\circ\psi$ belongs to $\F (Q,T)$ (if~and) only if
$\psi\in \F (R,T)\,.$\/}
\smallskip
\noindent
£2.1.2\quad {\it $\F_P (P)$ is a Sylow $p\-$subgroup of $\F (P)\,.$\/}
\smallskip
\noindent
Let us say that a subgroup $Q$ of $P$ is {\it fully centralized in $\F$\/} if for any $\F\-$morphism $\xi\,\colon Q\.C_P(Q)\to P$ we have $\xi \big(C_P (Q)\big) = C_P\big(\xi (Q)\big)\,.$ 
\smallskip
\noindent
£2.1.3\quad {\it For any  subgroup $Q$ of $P$ fully centralized in $\F\,,$ 
any $\F\-$morphism $\varphi\,\colon Q\to P$ and any subgroup $R$ of $N_P\big(\varphi(Q)\big)$ containing $\varphi (Q)$ such that $\F_P(Q)$ contains the action of $\F_R \big(\varphi(Q)\big)$ over $Q$ via $\varphi\,,$ there is an $\F\-$mor-phism 
$\zeta\,\colon R\to P$ fulfilling $\zeta\big(\varphi (u)\big) = u$ for any $u\in Q\,.$\/}

\medskip
£2.2. With the notation in~£1.1 above, it follows from [5,~Theorem~3.7] that $\F_{(b,G)}$
is a Frobenius $P\-$category. Moreover, we say that a subgroup $Q$ of $P$ is {\it $\F\-$nilcentralized\/}
 if, for  any $\varphi\in \F (P,Q)$ such that $Q' = \varphi (Q)$ is fully centralized in~$\F\,,$
the $C_P (Q')\-$categories $C_\F (Q')$ [5,~2.14] and $\F_{C_P (Q')}$ 
coincide; note that, according to [5,~Proposition~7.2], in $\F_{(b,G)}$ this definition agree with
£1.2 above. Similarly, we say that $Q$ is {\it $\F\-$selfcentralizing\/}  if we have
$$C_P\big (\varphi (Q))\i \varphi (Q)
\eqno £2.2.1\phantom{.}$$
 for any $\varphi \in \F (P,Q)\,;$ once again, according to [5,~Corollary~7.3], in $\F_{(b,G)}$ this definition agree with £1.2 above. Finally,  we say that a subgroup $R$ of $P$ is {\it $\F\-$radical\/} if it is
$\F\-$selfcentralizing and we have
$${\bf O}_p\big(\tilde\F (R)\big) = \{1\}
\eqno £2.2.2\phantom{.}$$
where   $\tilde\F (R) = \F (R)/\F_R (R)$ [5, 1.3].
  We respectively denote by $\F^{^{\rm nc}}\,,$  $\F^{^{\rm sc}}$ and $\F^{^{\rm rd}}$  the {\it full\/} subcategories of $\F$ over the respective sets of $\F\-$nilcentralized, $\F\-$selfcentralizing and $\F\-$radical subgroups of $P\,.$

\medskip 
£2.3.   We call {\it $\F^{^{\rm nc}}\!\-$chain\/} any functor  
$\frak q\,\colon \Delta_n\to \F^{^{\rm nc}}$ where the {\it $n\-$simplex\/} $\Delta_n$ is considered as a category with the morphisms --- denoted by~$i\bullet i'$ ---  defined by the order [5,~A2.2];  for any $\F\-$nilcentralized subgroup~$Q$ of~$P\,,$ let us denote by 
 $\frak q_Q\,\colon \Delta_0\to \F^{^{\rm nc}}$ the obvious 
 {\it $\F^{^{\rm nc}}\!\-$chain\/} sending $0$ to $Q\,.$ Following~[5,~A2.8],
 we denote  by $\ch^*(\F^{^{\rm nc}})$ the category  where the objects are all the 
$\F^{^{\rm nc}}\!\-$chains~$(\frak q,\Delta_n)$ and the morphisms from 
$\frak q\,\colon \Delta_n\to \F^{^{\rm nc}}$ to another $\F^{^{\rm nc}}\!\-$chain  $\frak r\,\colon \Delta_m\to \F^{^{\rm nc}}$ are the pairs $(\nu,\delta)$ 
formed by an {\it order preserving map\/} $\delta\,\colon \Delta_m\to \Delta_n$ and by a {\it natural isomorphism\/} $\nu\,\colon \frak q\circ\delta\cong \frak r\,,$   the composition being defined by the formula
$$(\mu,\varepsilon)\circ (\nu,\delta) = \big(\mu\circ (\nu *\varepsilon),\delta\circ \varepsilon\big)
\eqno £2.3.1.$$ 
Recall that we have a canonical functor [5, Proposition~A2.10]
$$\aut_{\F^{^{\rm nc}}} : \ch^*(\F^{^{\rm nc}})\too \Gr
\eqno £2.3.2\phantom{.}$$
mapping  any $\F^{^{\rm nc}}\!\-$chain $\frak q\,\colon \Delta_n\to \F^{^{\rm nc}}$ to the group of {\it natural automorphisms\/} of $\frak q\,.$ 

\medskip
£2.4. In  [6,~\S2] we introduce a {\it folded Frobenius $P\-$category\/} 
$(\F,\widehat\aut_{\F^{^{\rm sc}}})$ as a pair formed by a Frobenius $P\-$category $\F$ and  a functor 
$$\widehat\aut_{\F^{^{\rm sc}}} : \ch^*(\F^{^{\rm sc}})\too k^*\-\Gr
\eqno £2.4.1\phantom{.}$$
lifting the canonical functor  $\aut_{\F^{^{\rm sc}}}\,;$ here, we replace  {\it selfcentraling\/} by
{\it nilcentralized\/}: we call 
{\it folded Frobenius $P\-$category\/} $(\F,\widehat\aut_{\F^{^{\rm nc}}})$ a  pair formed by $\F$ and  a functor 
$$\widehat\aut_{\F^{^{\rm nc}}} : \ch^*(\F^{^{\rm nc}})\too k^*\-\Gr
\eqno £2.4.2\phantom{.}$$
lifting the canonical functor  $\aut_{\F^{^{\rm nc}}}\,;$ we also say that
$\widehat\aut_{\F^{^{\rm nc}}}$ is a {\it folder structure\/} of $\F\,.$ With the notation of~£1.1
above, Theorem~11.32 in [5] exhibits a {\it folder structure\/} of $\F_{(b,G)}\,,$ namely 
a functor $\widehat\aut_{(\F_{(b,G)})^{^{\rm nc}}}$ lifting $\aut_{(\F_{(b,G)})^{^{\rm nc}}}\,,$
that we call {\it Brauer folder structure\/} of $\F_{(b,G)}\,.$ Actually, both definitions coincide since any functor $\widehat\aut_{\F^{^{\rm sc}}}$ lifting $\aut_{\F^{^{\rm sc}}}$ can be extended to a unique functor $\widehat\aut_{\F^{^{\rm nc}}}$ lifting $\aut_{\F^{^{\rm nc}}}\,,$ as it shows our next result.

\bigskip
\noindent
{\bf Theorem~£2.5.} {\it Any functor $\widehat\aut_{\F^{^{\rm sc}}}$ lifting $\aut_{\F^{^{\rm sc}}}$ to the category $k^*\-\Gr$ can be extended to a unique functor lifting $\aut_{\F^{^{\rm nc}}}$
$$\widehat\aut_{\F^{^{\rm nc}}} : \ch^* (\F^{^{\rm nc}})\too k^*\-\Gr 
\eqno  £2.5.1.$$\/}
\par
\noindent
{\bf Proof:} Let $\frak X$ be a set of $\F\-$nilcentralized subgroups of $P$ which contains all the 
$\F\-$selfcentralizing subgroups of $P$ and is stable by $\F\-$isomorphisms; denoting by $\F^{^\frak X}$
the {\it full\/} subcategory of $\F$ over~$\frak X\,,$ assume that $\widehat\aut_{\F^{^{\rm sc}}}$ can be extended to a unique functor
$$\widehat\aut_{\F^{^{\frak X}}} : \ch^* (\F^{^{\frak X}})\too k^*\-\Gr 
\eqno £2.5.2.$$
\noindent
Assuming that $\frak X$ does not coincide with the set of all the $\F\-$nilcentralized~subgroups  
of~$P\,,$ let $V$ be a maximal  $\F\-$nilcentralized subgroup  which is not in~$\frak X\,;$  denoting by 
$\frak Y$ the union of $\frak X$ with all the subgroups of $P$  $\F\-$isomorphic to~$V\,,$ it is clear that it suffices to prove that $\widehat\aut_{\F^{^{\frak X}}}$ admits a unique extension 
to~$\ch^* (\F^{^{\frak Y}})\,.$

\smallskip
 For any chain  $\frak q\,\colon { \Delta}_n\to \F^{^{\frak Y}}\,,$  we choose an 
 $\F\-$morphism  $\alpha\,\colon \frak q (n)\to P$ such that $\alpha\big(\frak q (n)\big)$ is 
 {\it fully centralized\/}  in $\F$ [5,~Proposition~2.7] and denote\break by $\frak q^\alpha\,\colon \Delta_{n+1}\to \F^{^{\frak Y}}$ the chain which extends $\frak q$ and which maps  $n+1$ on $\alpha \big(\frak q (n)\big)\.C_P \Big(\alpha\big(\frak q (n)\big)\Big)$
and $(n\bullet n+1)$ on the $\F\-$morphism from $\frak q (n)$ to $\alpha \big(\frak q (n)\big)\.
C_P \Big(\alpha\big(\frak q (n)\big)\Big)$ induced by $\alpha\,;$ we have an obvious 
$\ch^* (\F^{^{\frak Y}})\-$morphism [5, A3.1]
$$({\rm id}_\frak q,\delta^n_{n+1}) : (\frak q^\alpha,\Delta_{n+1})\too (\frak q,\Delta_n)
\eqno £2.5.3\phantom{.}$$
 and the functor $\aut_{\F^{^{\frak Y}}}$ maps $({\rm id}_\frak q,\delta^n_{n+1})$ on a group homomorphism
$$\F(\frak q^\alpha)\too \F (\frak q)
\eqno £2.5.4\phantom{.}$$
which is surjective since  any $\sigma\in \F (\frak q)\i \F\big(\frak q (n)\big)$ can be
``extended'' to an $\F\-$automorphism of $\frak q^\alpha (n+1)$ [5,~statement~2.10.1].

\smallskip
Then, since $\alpha\big(\frak q (n)\big)$ is  is {\it $\F-$nilcentralized\/} and {\it fully centralized\/}  in $\F\,,$ the kernel of  homomorphism~£2.5.4 is a $p\-$group [5,~Corollary~4.7]; moreover, since $\frak q^\alpha (n+1)$ belongs to~$\frak X\,,$ the functor $\widehat\aut_{\F^{^\frak X}}$ and the structural inclusion 
$\F(\frak q^\alpha)\i \F\big(\frak q^\alpha (n+1)\big)$ determine a $k^*\-$subgroup 
$$\hat\F(\frak q^\alpha)\i \hat\F\big(\frak q^\alpha (n+1)\big) = \widehat\aut_{\F^{^\frak X}}
\big(\frak q^\alpha (n+1)\big)
\eqno £2.5.5\phantom{.}$$
 and, since  the kernel of  homomorphism~£2.5.4 is a $p\-$group,  this $k^*\-$subgroup induces  
a  central $k^*\-$extension  $\hat\F (\frak q)$ of~$\F (\frak q)$ such that we have a $k^*\-$group homomorphism
$$\hat\F(\frak q^\alpha)\too \hat\F (\frak q)
\eqno £2.5.6\phantom{.}$$
lifting homomorphism~£2.5.4. 
\eject

\smallskip
Note that, for a different choice $\alpha'\,\colon \frak q (n)\to P$ 
of~$\alpha\,,$ we have an $\F\-$iso-morphism $\alpha\big(\frak q (n)\big)\cong 
\alpha'\big(\frak q (n)\big)$ which can be extended to an $\F\-$isomorphism 
$\frak q^\alpha (n+1)\cong \frak q^{\alpha'} (n+1)$ [5,~statement~2.10.1]  and then $\widehat\aut_{\F^{^\frak X}}$ determines a $k^*\-$isomorphism 
$$\widehat\aut_{\F^{^\frak X}}\big(\frak q^\alpha (n+1)\big)\cong \widehat\aut_{\F^{^\frak X}}\big(\frak q^{\alpha'} (n+1)\big) 
\eqno £2.5.7\phantom{.}$$
 mapping 
$\hat\F(\frak q^{\alpha})$ onto $\hat\F(\frak q^{\alpha'})\,;$ moreover, it follows from [5,~Proposition~4.6] that two such $\F\-$isomorphisms are $C_P \Big(\alpha'\big(\frak q (n)\big)
\Big)\-$conjugate and therefore our definition of $\hat\F (\frak q)$ does not depend on our choice of 
$\alpha\,.$ Similarly, if $\frak q (n)$ belongs to $\frak X$ then the functor $\widehat\aut_{\F^{^\frak X}}$ already defines a $k^*\-$group $\widehat\aut_{\F^{^\frak X}}\big(\frak q (n)\big)$ 
and, denoting by $\frak q^\alpha_{n,n+1}\,\colon \Delta_1\to \F^{^\frak X}$ the chain mapping $0$
on $\frak q (n)\,,$ $1$ on $\frak q^\alpha (n+1)$ and $(0\bullet 1)$ on 
$\frak q^\alpha (n\bullet n+1)\,,$ also defines a $k^*-$group homomorphism
$$\widehat\aut_{\F^{^\frak X}}(\frak q^\alpha_{n,n+1})_{}\too \widehat\aut_{\F^{^\frak X}}\big(\frak q (n)\big)
\eqno £2.5.8$$
inducing a {\it canonical\/} $k^*\-$group isomorphism from $\hat\F (\frak q)$ in~£2.5.6 above onto the inverse~image
of $\aut_{\F^\frak Y} (\frak q)\i \aut_{\F^\frak X}\big(\frak q (n)\big)$ in $\widehat\aut_{\F^{^\frak X}}\big(\frak q (n)\big)\,;$ in particular, if the image of $\frak q$ is contained in $\frak X\,,$  we get a {\it canonical\/} $k^*\-$group isomorphism $\hat\F (\frak q) \cong \widehat\aut_{\F^{^\frak X}}(\frak q)\,.$

\smallskip
 Now, for any $\ch^*(\F^{^\frak Y})\-$morphism $(\nu,\delta)\,\colon (\frak r,\Delta_m)
 \to (\frak q,\Delta_n)\,,$ choosing suitable $\F\-$morphisms  $\alpha\,\colon \frak q (n)\to P$
 and $\beta\,\colon \frak r (m)\to P$ as above, we have to exhibit a $k^*\-$group homomorphism 
 $\hat\F (\frak r)\to \hat\F (\frak q)$ lifting $\aut_{\F^{^\frak Y}} (\nu,\delta)\,.$ 
 Firstly, we assume that the image of $\frak r(\delta (n))$ {\it via\/} 
 $\frak r (\delta (n)\bullet m)$ is {\it normal\/} in $\frak r (m)\,;$ in this case,  
 $\beta\Big(\frak r (\delta (n)\bullet m)\big(\frak r(\delta (n))\big)\Big)$ is normal in 
 $\frak r^\beta (m+\!1)$ and, according to~[5,~statement~2.10.1], there is an  
 $\F\-$morphism 
 $$\hat\nu :  \frak r^\beta (m+\!1) \too N_P\Big(\alpha\big(\frak q (n)\big)\Big)
 \eqno~£2.5.9\phantom{.}$$
  extending  the  $\F\-$morphism 
$$\beta\Big(\frak r (\delta (n)\bullet m)\big(\frak r(\delta (n))\big)\Big)\cong \frak r(\delta (n))
\buildrel \nu_n\over\cong \frak q(n) \cong \alpha\big(\frak q (n)\big)\i P
\eqno~£2.5.10,$$
and we set $U = \hat\nu \big(\frak r^\beta (m+\!1)\big)\.C_P\Big(\alpha\big(\frak q (n)\big)\Big)\,.$
  Then, we consider the chains
$$\frak q^{\alpha,\nu} : \Delta_{n+2}\too \F^{^\frak Y}\qq 
\frak r^{\beta,\nu} : \Delta_{m+2}\too \F^{^\frak Y}
\eqno £2.5.11\phantom{.}$$
respectively extending the chains $\frak q^\alpha$ and $\frak r^\beta$ defined above, fulfilling 
$$\frak q^{\alpha,\nu}(n+2) = U = \frak r^{\beta,\nu}(m+2)
\eqno £2.5.12\phantom{.}$$ 
and, since $\alpha\big(\frak q (n)\big)\i\hat\nu \Big(\beta\big(\frak r (m)\big)\Big)\,,$ respectively mapping $(n\!+\!1\bullet n\!+\!2)$ and $(m\!+\!1\bullet m\!+\!2)$ on
the inclusion $\frak q^\alpha (n+1) \i U$ and on the $\F\-$morphism from\break 
\eject
\noindent
$\frak r^\beta (m+1)$~to $U$ induced by $\hat\nu\,$ . Note that, since the centralizer of 
$\alpha\big(\frak q (n)\big)$ contains  $C_P \bigg(\hat\nu \Big(\beta\big(\frak r (m)\big)\Big)\bigg)$ and $\beta\big(\frak r(m)\big)$ is fully centralized in $\F\,,$ we still have  $U = \hat\nu\big(\beta (\frak r (m))\big)\. C_P\Big(\alpha\big(\frak q (n)\big)\Big)\,.$
Moreover, it follows from [5,~Proposition~4.6] that another choice 
$\hat\nu'$ of the $\F\-$morphism~£2.5.9 is $C_P\Big(\alpha\big(\frak q (n)\big)\Big)\-$conjugate 
of $\hat\nu$ and, in particular, the group $U$ does not change.

\smallskip
With all this notation, we have obvious $\ch^*(\F^{^\frak Y})\-$morphisms
$$\eqalign{({\rm id}_{\frak q^\alpha},\delta^{n+1}_{n+2}) : (\frak q^{\alpha,\nu},\Delta_{n+2})&\too (\frak q^\alpha,\Delta_{n+1}) \cr 
({\rm id}_{\frak r^\beta},\delta^{m+1}_{m+2}) :(\frak r^{\beta,\nu},\Delta_{m+2})&\too (\frak r^\beta,\Delta_{m+1})\cr}
\eqno £2.5.13\phantom{.}$$
and, considering  the maps
$$\Delta_{n+2}\buildrel \sigma_n\over\longleftarrow\Delta_1\buildrel \sigma_m\over\too \Delta_{m+2}\qq \Delta_{n+1}\buildrel \tau_n\over
\longleftarrow\Delta_0\buildrel \tau_m\over\too \Delta_{m+1}
\eqno £2.5.14\phantom{.}$$
respectively mapping $i$ on $i +n+1$ and $i +m+1\,,$ the  
 $\ch^*(\F^{^\frak Y})\-$morphisms above  determine the following 
 $\ch^*(\F^{^\frak X})\-$morphisms
$$(\frak q^{\alpha,\nu}\!\circ \sigma_n,\Delta_{1})\too (\frak q^\alpha\circ \tau_n,
\Delta_{0}) \!\!\qq\!\!  (\frak r^{\beta,\nu}\!\circ \sigma_m,\Delta_{1})\too (\frak r^\beta \circ
\tau_m,\Delta_{0})
\eqno £2.5.15.$$
Then, the functor $\widehat\aut_{\F^{^\frak X}}$ maps these morphisms on 
$k^*\-$group homomorphisms
$$\hat\F (\frak q^{\alpha,\nu}\!\circ \sigma_n)\too \hat\F (\frak q^\alpha \circ  \tau_n) \qq 
\hat\F (\frak r^{\beta,\nu}\!\circ \sigma_m)\too \hat\F (\frak r^\beta \circ \tau_m) 
\eqno £2.5.16.$$

\smallskip
But note that $\F(\frak q^{\alpha,\nu})\,,$ $\F(\frak q^\alpha)\,,$   $\F(\frak r^{\beta,\nu})$ 
and $\F ( r^\beta)$ are respectively contained in $\F (\frak q^{\alpha,\nu}\!\circ \sigma_n)\,,$ 
$\F (\frak q^\alpha \circ \tau_n)\,,$ $\F (\frak r^{\beta,\nu}\!\circ \sigma_m)$
and $\F (\frak r^\beta \circ \tau_m) \,,$ and therefore,  considering the corresponding
inverse images in $\hat\F (\frak q^{\alpha,\nu}\!\circ \sigma_n)\,,$ 
$\hat\F (\frak q^\alpha \circ \tau_n)\,,$ $\hat\F (\frak r^{\beta,\nu}\!\circ \sigma_m)$
and $\hat\F (\frak r^\beta \circ \tau_m) \,,$ the $k^*\-$group homomorphisms~£2.5.16 
induce $k^*\-$group homomorphisms (cf.~£2.5.8)
$$\hat\F (\frak q^{\alpha,\nu})\too \hat\F (\frak q^\alpha ) \qq \hat\F (\frak r^{\beta,\nu}) \too 
\hat\F (\frak r^\beta) 
\eqno £2.5.17.$$
More explicitly, we actually have
$$\F (\frak q^{\alpha,\nu} \circ  \sigma_n)  = \F (U) = \F (\frak q^{\beta,\nu} \circ \sigma_m) 
\eqno £2.5.18\phantom{.}$$
and the structural inclusions $\F (\frak q^{\alpha,\nu})\i \F(U)$ and $\F (\frak r^{\beta,\nu})\i \F(U)$ 
induce an inclusion $\F (\frak r^{\beta,\nu})\i \F (\frak q^{\alpha,\nu})\,;$ indeed, an element $\theta$
in $\F (\frak r^{\beta,\nu})$ stabilizes the subgroups $\hat\nu\Big(\beta\big(\frak r (i\bullet m)\big(\frak r(i)\big)\big)\Big)$ of $U$
for any $i\in \Delta_m\,,$ so that it stabilizes 
$$\alpha \big(\frak q (n)\big) = \hat\nu\bigg(\beta\Big(\frak r \big(\delta (n)\bullet m\big)
\big(\frak r(\delta (n))\big)\Big)\bigg)
\eqno £2.5.19,$$
and therefore $\theta$ also stabilizes $C_P\Big(\alpha \big(\frak q (n)\big)\Big) = C_U\Big(\alpha \big(\frak q (n)\big)\Big)\,;$ thus, it stabilizes the subgroup $\frak q^\alpha (n+1)$ of $U$ and therefore
$\theta$ belongs to $\F (\frak q^{\alpha,\nu})\,.$
\eject

\smallskip
Moreover, we claim that 
$$\big(\aut_{\F^\frak Y} ({\rm id}_{\frak r^\beta},\delta^{m+1}_{m+2})\big)
\big(\F (\frak r^{\beta,\nu})\big) = \F (\frak r^\beta)   
\eqno £2.5.20.$$
Indeed, an element $\theta$ in $\F (\frak r^\beta)$ acts on $\beta \big(\frak r (m)\big)$ 
determining an automorphism $\hat\theta$ of $\hat\nu\Big(\beta \big(\frak r (m)\big)\Big)$
and, as above, this automorphism stabilizes $\alpha \big(\frak q (n)\big)$ inducing an $\F\-$morphism
$$\eta : \alpha \big(\frak q (n)\big)\cong \alpha \big(\frak q (n)\big)\i P
\eqno £2.5.21;$$
but, we are assuming that $\alpha \big(\frak q (n)\big)$ is normal in $\hat\nu\Big(\beta \big(\frak r (m)\big)\Big)\,,$ so that this group is normal in $\frak r^{\beta,\nu} (m+2)$ (cf.~£2.5.12). Hence, it follows
from [5,~statement~2.10.1] that $\eta$ can be extended to an $\F\-$morphism 
$\hat\eta\,\colon \frak r^{\beta,\nu} (m+2)\to P\,;$ then, the restriction of $\hat\eta$ to
$\hat\nu\Big(\beta \big(\frak r (m)\big)\Big)$ and the $\F\-$morphism 
$$\hat\nu\Big(\beta \big(\frak r (m)\big)\Big)\buildrel \hat\theta\over\cong \hat\nu\Big(\beta \big(\frak r (m)\big)\Big)\i P
\eqno £2.5.22\phantom{.}$$
coincide over the subgroup $\alpha \big(\frak q (n)\big)$ and therefore, according to [5,~Proposition~4.6],  these homomorphisms are $C_P \Big(\alpha \big(\frak q (n)\Big)\-$conjugate.
In conclusion, up to a modification in our choice of $\hat\eta\,,$ we may assume that the restriction of  
$\hat\eta$ to $\hat\nu\Big(\beta \big(\frak r (m)\big)\Big)$ coincides with $\hat\theta$ and therefore
that $\hat\eta$ stabilizes $\hat\nu\big( \frak r^{\beta,\nu} (m+1)\big)$ and 
$\hat\nu\big( \frak r^{\beta,\nu} (m+2)\big)\,,$ so that $\hat\eta$ induces an element
of $\F (\frak r^{\beta,\nu})$ lifting $\theta\,.$

\smallskip
Consequently, we have the following commutative diagram
$$\matrix{\F (U)&\j &\F (\frak q^{\alpha,\nu})&\too &\F (\frak q^\alpha)&\too &\F (\frak q)\cr
\Vert&\phantom{\big\uparrow}& \cup&&&&\hskip-40pt{\scriptstyle \aut_{\F^{^\frak Y}} (\nu,\delta)}\uparrow\cr
\F (U)&\j &\F (\frak r^{\beta,\nu})&\too &\F (\frak r^\beta)&\too &\F(\frak r)\cr}
\eqno £2.5.23;$$
Moreover, since $\frak q^\alpha (n+1)$ and $\frak r^\beta (m+1)$ are $\F\-$selfcentralizing,
the kernels of the compositions of the horizontal arrows  are 
$\F_{C_U (\alpha (\frak q (n)))} (U)$ for the top and $\F_{C_U (\hat\nu (\beta (\frak r (m))))} (U)$
for the bottom, and the bottom composition is surjective; hence, since 
$\F_{C_U (\hat\nu (\beta (\frak r (m))))} (U)$ is contained in $\F_{C_U (\alpha (\frak q (n)))} (U)$
and they respectively lift canonically to $\hat\F (\frak r^{\beta,\nu})$ and to 
$\hat\F (\frak q^{\alpha,\nu})$ [5, Corollaire 4.7], we get a {\it unique\/} $k^*\-$group homomorphism
$$\widehat{\aut}_{\F^{^\frak Y}} (\nu,\delta) : \hat\F (\frak r)\too \hat\F (\frak q)
\eqno £2.5.24\phantom{.}$$
lifting $\aut_{\F^{^\frak Y}} (\nu,\delta)$ and such that the corresponding diagram of $k^*\-$group
homomorphisms
$$\matrix{\hat\F (U)&\j &\hat\F (\frak q^{\alpha,\nu})&\too &\hat\F (\frak q^\alpha)&\too &\hat\F (\frak q)\cr
\Vert&\phantom{\big\uparrow}& \cup&&&
&\hskip-40pt{\scriptstyle \widehat{\aut}_{\F^{^\frak Y}} (\nu,\delta)}\uparrow\cr
\hat\F (U)&\j &\hat\F (\frak r^{\beta,\nu})&\too &\hat\F (\frak r^\beta)&\too &\hat\F(\frak r)\cr}
\eqno £2.5.25\phantom{.}$$
is commutative.
\eject

\smallskip
 Consider another $\ch^*(\F^{^\frak Y})\-$morphism $(\mu,\varepsilon)\,\colon 
 (\frak t,\Delta_\ell)\to (\frak r,\Delta_m)\,,$ so that 
 $$(\nu,\delta) \circ (\mu,\varepsilon) = \big(\nu\circ (\mu * \delta), \varepsilon\circ\delta\big)
 \eqno £2.5.26\phantom{.}$$
and set $\lambda = \nu\circ (\mu * \delta)$ and $\varphi = \varepsilon\circ\delta\,;$
then, choosing a suitable $\F\-$morphism $\gamma\,\colon \frak t (\ell)\to P$ as above,
we still assume  that  the images of $\frak t(\varphi (n))$ {\it via\/} 
 $\frak t (\varphi (n)\bullet \ell)$ and of $\frak t(\varepsilon (m))$ {\it via\/} 
 $\frak t (\varepsilon (m)\bullet \ell)$ are {\it normal\/} in $\frak t (\ell)\,.$ In particular, 
 this implies that  the image of $\frak r(\delta (n))$ {\it via\/} 
 $\frak r (\delta (n)\bullet m)$ is {\it normal\/} in $\frak r (m)\,;$ that is to say, we  have already defined 
 the $k^*\-$group homomorphisms $\widehat{\aut}_{\F^{^\frak Y}} (\nu,\delta)\,,$
 $\widehat{\aut}_{\F^{^\frak Y}} (\mu,\varepsilon)$ and 
 $\widehat{\aut}_{\F^{^\frak Y}} (\lambda,\varphi) $ respectively lifting 
 $\aut_{\F^{^\frak Y}} (\nu,\delta)\,,$ $\aut_{\F^{^\frak Y}} (\mu,\varepsilon)$ and 
 $\aut_{\F^{^\frak Y}} (\lambda,\varphi) $ and we want to prove that
 $$\widehat{\aut}_{\F^{^\frak Y}} (\lambda,\varphi) = \widehat{\aut}_{\F^{^\frak Y}} (\nu,\delta)\circ \widehat{\aut}_{\F^{^\frak Y}} (\mu,\varepsilon)
 \eqno £2.5.27.$$
 
 \smallskip
More explicitly, applying the construction in~£2.5.9 above to the $\ch^*(\F^{^\frak Y})\-$ morphisms
$(\nu,\delta)\,,$ $(\mu,\varepsilon)$ and $(\varphi,\lambda)\,,$ we get  $\F\-$morphisms 
$$\eqalign{\hat\nu :  \frak r^\beta (m+\!1)& \too N_P\Big(\alpha\big(\frak q (n)\big)\Big)\cr
\hat\mu : \frak t^\gamma (\ell+1)&\too N_P \Big(\beta \big(\frak r (m)\big)\Big)\cr
\hat\lambda : \frak t^\gamma (\ell+1)&\too N_P \Big(\alpha \big(\frak q (n)\big)\Big)\cr}
\eqno £2.5.28;$$
actually, it is clear that the respective images of $\hat\nu\,,$ $\hat\mu$ and $\hat\lambda$ are 
respectively contained in $\frak q^\alpha (n+1)\,,$ $\frak r^\beta (m+1)$ and $\frak q^\alpha (n+1)$
and, with evident notation, our construction can be explicited in the following commutative diagram
$$\matrix{\frak t (\ell)\hskip-6pt&\cong &\hskip-6pt\gamma \big(\frak t (\ell)\big)&\hskip-6pt\i &\hskip-6pt\frak t^\gamma (\ell+1) &&\hskip-6pt\buildrel \hat\lambda\over\too 
&&\hskip-6pt\frak q^\alpha (n+1)\cr
&&&&\hskip-6pt\Vert\cr
\big\uparrow\hskip-6pt&&&&\hskip-6pt\frak t^\gamma  (\ell+1)&\hskip-6pt\buildrel \hat\mu\over\too  &\hskip-6pt\frak r^\beta (m+1)\cr 
&&&&&&\hskip-6pt\Vert\cr
\frak t(\varepsilon (m))\hskip-6pt&\buildrel \mu_m\over\cong \hskip-6pt&\hskip-6pt\frak r (m)&\hskip-6pt\cong &\hskip-6pt\beta\big(\frak r (m)\big)&\hskip-6pt\i &\hskip-6pt\frak r^\beta (m+1)&\hskip-6pt\buildrel \hat\nu\over\too &\hskip-6pt\frak q^\alpha (n+1)\cr
\big\uparrow&&\hskip-6pt\big\uparrow&&&&&&\hskip-6pt\cup\cr
\frak t \big(\varphi (n)\big)\hskip-6pt&\buildrel \mu_{\delta (n)}\over \cong
&\hskip-6pt\frak r\big(\delta (n)\big) &\hskip-6pt\buildrel \nu_n\over\cong &\hskip-6pt\frak q (n) 
&&\hskip-6pt\cong &&\hskip-6pt\alpha\big(\frak q (n)\big)\cr}
\eqno £2.5.29.$$

\smallskip
That is to say, according to~£2.5.10 above, $\hat\lambda\,,$  $\hat\mu$ and $\hat\nu$ respectively extend the 
$\F\-$morphisms
$$\eqalign{\gamma\Big(\frak t (\varphi (n)\bullet \ell)\big(\frak t(\varphi (n))\big)\Big)
&\cong \frak t(\varphi (n))\buildrel \lambda_n\over\cong \frak q(n) \cong \alpha\big(\frak q (n)\big)\i P\cr
\gamma\Big(\frak t (\varepsilon (m)\bullet \ell)\big(\frak t(\varepsilon (m))\big)\Big)
&\cong \frak t(\varepsilon (m))\buildrel \mu_m\over\cong \frak r (m) \cong \beta\big(\frak r (m)\big)\i P\cr
\beta\Big(\frak r (\delta (n)\bullet m)\big(\frak r(\delta (n))\big)\Big)&\cong \frak r(\delta (n))
\buildrel \nu_n\over\cong \frak q(n) \cong \alpha\big(\frak q (n)\big)\i P\cr}
\eqno~£2.5.30\phantom{.}$$
and, since $\beta\Big(\frak r (\delta (n)\bullet m)\big(\frak r(\delta (n))\big)\Big)$ is contained in
$\beta\big(\frak r (m)\big)\,,$ it is easily checked that the composition $\hat\nu\circ \hat\mu$
also extends the top $\F\-$morphism in~£2.5.30; then, as above, it follows from [5,~Proposition~4.6] that
$\hat\lambda$ and $\hat\nu\circ \hat\mu$ are $C_P\Big(\alpha\big(\frak q (n)\big)\Big)\-$conjugate;
actually, up to a modification of our choice of $\hat\lambda\,,$ we may assume that they coincide.

 \smallskip
Moreover, we have to consider  chains
$$\eqalign{\frak q^{\alpha,\nu,\lambda} &: \Delta_{n+3}\too \F^{^\frak Y}\cr 
\frak r^{\beta,\mu,\nu} &: \Delta_{m+3}\too \F^{^\frak Y}\cr
 \frak t^{\gamma,\mu,\nu} &: \Delta_{\ell+3}\too \F^{^\frak Y}\cr}
\eqno £2.5.31\phantom{.}$$
respectively extending the chains $\frak q^{\alpha,\nu}\,,$ $\frak r^{\beta,\mu}$ and 
$ \frak t^{\gamma,\mu}\,;$ recall that (cf.~£2.5.12)
$$\eqalign{\frak q^{\alpha,\nu} (n+2) &= \hat\nu \Big(\beta \big(\frak r (m)\big)\Big)
\.C_P\Big(\alpha\big(\frak q (n)\big)\Big)\cr 
\frak r^{\beta,\mu}(m+ 2) &= \hat\mu\Big( \gamma\big(\frak t (\ell )\big)\Big)
\. C_P\Big(\beta \big(\frak r (m)\big)\Big)= \frak t^{\gamma,\mu} (\ell+ 2)\cr}
\eqno £2.5.32\phantom{.}$$
and that, according to our remark above and since we assume that 
$\hat\lambda = \hat\nu\circ\hat\mu\,,$ we still have
$$\frak q^{\alpha,\lambda} (n +2) = \hat\nu\bigg(\hat\mu \Big(\gamma\big(\frak t (\ell )\big)\Big)\bigg)\. C_P\Big(\alpha\big(\frak q (n)\big)\Big)
\eqno £2.5.33;$$
thus, since $\beta \big(\frak r (m)\big)\i \hat\mu \Big(\gamma\big(\frak t (\ell )\big)\Big)\,,$
we get $\frak q^{\alpha,\nu} (n+2)\i \frak q^{\alpha,\lambda} (n +2)$ and, since 
the centralizer of $\alpha\big(\frak q (n)\big)$ contains the centralizer of $\hat\nu \Big(\beta\big(\frak r (m)\big)\Big)\,,$  $\hat\nu$ induces an $\F\-$morphism
$$\frak r^{\beta,\mu}(m+ 2)  = \frak t^{\gamma,\mu} (\ell+ 2)\too \frak q^{\alpha,\lambda} (n +2) 
\eqno £2.5.34;$$
then, we complete our definition of $\frak q^{\alpha,\nu,\lambda} \,,$ $\frak r^{\beta,\mu,\nu}$
and $\frak t^{\gamma,\mu,\nu}$ by setting
$$\frak q^{\alpha,\nu,\lambda} (n+3) = \frak r^{\beta,\mu,\nu} (m+3) 
= \frak t^{\gamma,\mu,\nu}(\ell +3) = \frak q^{\alpha,\lambda} (n +2)
\eqno £2.5.35,$$ 
and respectively mapping  $(n +2\bullet n+3)\,,$  
$(m + 2\bullet m+3)$ and $(\ell + 2\bullet \ell + 3)$ on the inclusion 
$\frak q^{\alpha,\nu} (n + 2)\i \frak q^{\alpha,\lambda} (n +2)$ and on the $\F\-$morphism~£2.5.34
induced by $\hat\nu\,.$

\smallskip
Now, it is clear that the functor $\aut_{\F^\frak Y}$ applied to the obvious 
$\ch^*(\F^{^\frak Y})\-$ morphisms
$$\eqalign{({\rm id}_{\frak q^{\alpha,\nu}},\delta^{n+2}_{n+3}) : 
(\frak q^{\alpha,\nu,\lambda},\Delta_{n+3})&\too (\frak q^{\alpha,\nu},\Delta_{n+2}) \cr 
({\rm id}_{\frak r^{\beta,\mu}},\delta^{m+2}_{m+3}) :(\frak r^{\beta,\mu,\nu},\Delta_{m+3})
&\too (\frak r^{\beta,\mu},\Delta_{m+2})\cr
({\rm id}_{\frak t^{\gamma,\mu}},\delta^{m+2}_{m+3}) :(\frak t^{\gamma,\mu,\nu},\Delta_{m+3})
&\too (\frak t^{\gamma,\mu},\Delta_{m+2})\cr}
\eqno £2.5.36\phantom{.}$$
yields group homomorphisms
$$\F (\frak q^{\alpha,\nu,\lambda})\to \F (\frak q^{\alpha,\nu})\quad \!\!\!\!,\!\!\!\!\quad
\F (\frak r^{\beta,\mu,\nu})\to \F (\frak r^{\beta,\mu})\quad\!\!\!\!,\!\!\!\!\quad 
\F (\frak t^{\gamma,\mu,\nu})\to \F (\frak t^{\gamma,\mu})
\eqno £2.5.37;$$
as in~£2.5.16 above, considering the maps
$$\eqalign{\hat\sigma_n : \Delta_1\too \Delta_{n+3}&\qq \hat\tau_n : \Delta_0\too \Delta_{n+2}\cr
\hat\sigma_m : \Delta_1\too \Delta_{m+3}&\qq \hat\tau_m : \Delta_0\too \Delta_{m+2}\cr
\hat\sigma_\ell : \Delta_1\too \Delta_{\ell+3}&\qq \hat\tau_\ell : \Delta_0\too \Delta_{\ell+2}\cr}
\eqno £2.5.38\phantom{.}$$
respectively sending $i$ to $i + n + 2\,,$ to $i + m + 2$ and to $i + \ell + 2\,,$  the functor $\widehat\aut_{\F^{^\frak X}}$ still induces $k^*\-$group homomorphisms
$$\eqalign{\hat\F (\frak q^{\alpha,\nu,\lambda})&\too \hat\F (\frak q^{\alpha,\nu})\cr
\hat\F (\frak r^{\beta,\mu,\nu})&\too \hat\F (\frak r^{\beta,\mu})\cr
\hat\F (\frak t^{\gamma,\mu,\nu})&\too \hat\F (\frak t^{\gamma,\mu})\cr}
\eqno £2.5.39;$$
moreover it is quite clear that $\hat\F (\hat\frak t^{\gamma,\mu,\nu}) = \hat\F (\hat\frak t^{\gamma,\nu})\,.$
Consequently, the functoriality of $\widehat\aut_{\F^{^\frak X}}$ guarantees the commutativity of
the following diagram
$$\matrix{\hat\F (\hat\frak t)&\hskip-5pt\leftarrow&\hskip-5pt
\hat\F (\hat\frak t^{\gamma,\mu,\nu})&\hskip-5pt  = &\hskip-5pt\hat\F (\hat\frak t^{\gamma,\nu})&\hskip-5pt&\hskip-5pt\i &\hskip-5pt&\hskip-5pt\hat\F (\hat\frak q^{\gamma,\lambda})&\hskip-10pt = &\hskip-8pt\hat\F (\hat\frak q^{\gamma,\lambda})\cr
\Vert&\hskip-5pt&\hskip-5pt\downarrow&\hskip-5pt\phantom{\big\downarrow}&\hskip-5pt&\hskip-5pt&\hskip-5pt\phantom{\Big\downarrow}&\hskip-5pt&\hskip-5pt\cup\cr
\hat\F (\hat\frak t)&\hskip-5pt\leftarrow&\hskip-5pt\hat\F (\hat\frak t^{\gamma,\mu})&\hskip-5pt\i
&\hskip-5pt\hat\F (\hat\frak r^{\beta,\lambda})&\hskip-5pt\leftarrow&\hskip-5pt\hat\F 
(\hat\frak r^{\beta,\mu,\nu})&\hskip-5pt\i
&\hskip-5pt\hat\F (\hat\frak q^{\alpha,\nu,\lambda})\cr
\downarrow&\hskip-5pt&\hskip-5pt&\hskip-5pt\phantom{\Big\downarrow}&\hskip-5pt\downarrow&\hskip-5pt&\hskip-5pt\downarrow&\hskip-5pt&\hskip-5pt\downarrow&\hskip-5pt&\hskip-5pt\downarrow\cr
\hat\F (\frak t)&\hskip-5pt&\hskip-5pt&\hskip-5pt&\hskip-5pt \hat\F (\hat\frak r)&\hskip-5pt\leftarrow&\hskip-5pt
\hat\F (\hat\frak r^{\beta,\mu})&\hskip-5pt\i  &\hskip-5pt\hat\F (\hat\frak q^{\alpha,\nu})\cr 
&\hskip-5pt&\hskip-5pt\searrow&\hskip-5pt&\hskip-5pt\downarrow&\hskip-5pt\phantom{\Bigg\downarrow}&\hskip-5pt&\hskip-5pt &\hskip-5pt\downarrow\cr
 &\hskip-5pt&\hskip-5pt&\hskip-5pt&\hskip-5pt\hat\F (\frak r)&\hskip-5pt&\hskip-5pt&\hskip-5pt&\hskip-5pt\hat\F (\hat\frak q)&\hskip-10pt =&\hskip-8pt\hat\F (\hat\frak q)\cr 
&\hskip-5pt&\hskip-5pt&\hskip-5pt&\hskip-5pt\phantom{\Bigg\downarrow}&\hskip-5pt&\hskip-5pt \searrow&\hskip-5pt&\hskip-5pt\downarrow\cr
&\hskip-5pt&\hskip-5pt&\hskip-5pt &\hskip-5pt&\hskip-5pt&\hskip-5pt&\hskip-5pt&\hskip-5pt\hat\F (\frak q)\cr}
\eqno £2.5.40;$$
thus, by uniqueness, in this case we obtain 
$$\widehat\aut_{\F^{^\frak Y}}(\nu,\delta)\circ\widehat\aut_{\F^{^\frak Y}}(\mu,\varepsilon) =\widehat\aut_{\F^{^\frak Y}}\big((\nu,\delta)\circ (\mu,\varepsilon)\big)
\eqno £2.5.41.$$

\smallskip
 Secondly, assume that the image of $\frak r\big(\delta (n)\big)$ by 
 $\frak r (\delta (n)\!\bullet\! m)$ is not normal in $\frak r (m)\,;$ let $m'$ be the maximal element in $\Delta_m - \Delta_{\delta (n) -1}$ such that  the image of 
 $\frak r\big(\delta (n)\big)$  by $\frak r (\delta (n)\!\bullet\! m')$ is normal in~$\frak r (m')$ and denote by $R_{(\nu,\delta)}$ the normalizer of the image of $\frak r\big(\delta (n)\big)$
 in $\frak r (m'+1)\,,$ by $\frak r_{(\nu,\delta)}
\colon \Delta_{m+1}\to \F^{^\frak Y}$ the functor fulfilling 
$$\frak r_{(\nu,\delta)}\circ \delta^m_{m'+1} = \frak r\qq
\frak r_{(\nu,\delta)}(m'+1) = R_{(\nu,\delta)}
\eqno £2.5.42\phantom{.}$$
 and mapping $(m'\!+\!1\bullet m'\!+\!2)$ on the inclusion map $R_{(\nu,\delta)}\to\frak r (m'+1)\,,$ and by $\frak r'_{(\nu,\delta)}$ the restriction of $\frak r_{(\nu,\delta)}$
to $\Delta_{m'+1}\,;$ then, it is quite clear that $\F (\frak r_{(\nu,\delta)}) 
= \F (\frak r)$ and it is easily checked that $\hat\F (\frak r_{(\nu,\delta)}) 
= \hat\F (\frak r)\,;$ moreover, we have an evident $\ch^*(\F^{^\frak Y})\-$morphism  
$$(\nu',\delta') : (\frak r'_{(\nu,\delta)},\Delta_{m'+1})\too (\frak q,
\Delta_{n})
\eqno £2.5.43\phantom{.}$$
such that 
$$(\nu',\delta')\circ ({\rm id}_{\frak r'_{(\nu,\delta)}},\iota^m_{m'}) =
(\nu,\delta)\circ  ({\rm id}_\frak r,\delta^m_{m'+1})
\eqno £2.5.44\phantom{.}$$
 where $\iota^m_{m'}\colon\Delta_{m'+1} \to \Delta_{m+1}$ denotes the natural inclusion,  we clearly  have\break $\widehat\aut_{\F^{^\frak Y}}({\rm id}_\frak r,\delta^m_{m'+1}) = 
{\rm id}_{\hat\F (\frak r)}$ and in £2.5.24 above we have already defined 
  $\widehat\aut_{\F^{^\frak Y}} (\nu',\delta')\,;$ on the other hand, arguing by induction on $\vert\frak r (m)\vert/\vert\frak q (n)\vert$, we may assume that $\widehat\aut_{\F^{^\frak Y}}
({\rm id}_{\frak r'_{(\nu,\delta)}},\iota^m_{m'})$ is already defined and then we set
$$\widehat\aut_{\F^{^\frak Y}}(\nu,\delta) =
\widehat\aut_{\F^{^\frak Y}} (\nu',\delta') \circ\widehat\aut_{\F^{^\frak Y}}
({\rm id}_{\frak r'_{(\nu,\delta)}},\iota^m_{m'})
\eqno £2.5.45.$$

\smallskip
  For another $\ch^*(\F^{^\frak Y})\-$morphism $(\mu,\varepsilon)\,\colon (\frak t,\Delta_\ell)\to (\frak r, \Delta_m)\,,$  we claim that
$$\widehat\aut_{\F^{^\frak Y}}(\nu,\delta)\circ\widehat\aut_{\F^{^\frak Y}}(\mu,\varepsilon) = \widehat\aut_{\F^{^\frak Y}}\big((\nu,\delta)\circ (\mu,\varepsilon)\big)
\eqno £2.5.46;$$
we argue by induction firstly on $\vert\frak t (\ell)\vert/\vert\frak q (n)\vert$ and after
on $\vert\frak t (\ell)\vert/\vert\frak r (m)\vert\,.$ First of all, we assume that the image of $\frak r\big(\delta (n)\big)$ in~$\frak r (m)$ by $\frak r\big(\delta (n)\bullet m\big)$ is not normal; with the notation above, denote by $\ell'$ the maximal element in 
$\Delta_\ell - \Delta_{(\varepsilon\circ\delta) (n) -1}$ such that  the image of $\frak t\big((\varepsilon\circ\delta) (n)\big)$ 
by  $\frak t \big((\varepsilon\circ\delta) (n)\bullet \ell'\big)$ is normal in
$\frak t (\ell')\,;$ then, it is clear that
$\varepsilon (m')\le \ell' < \varepsilon (m)$ and easily checked that we have a $\ch^*(\F^{^\frak Y})\-$morphism
$$(\mu_{(\nu,\delta)},\varepsilon_{(\nu,\delta)}) : (\frak t_{(\nu,\delta) \circ
(\mu,\varepsilon)},\Delta_{\ell +1})\too (\frak r_{(\nu,\delta)} ,\Delta_{m +1})
\eqno £2.5.47\phantom{.}$$
such that 
$$({\rm id}_\frak r,\delta^m_{m'+1})\circ (\mu_{(\nu,\delta)}, \varepsilon_{(\nu,\delta)}) = (\mu,\varepsilon)\circ ({\rm id}_\frak t,\delta^\ell_{\ell'+1})
\eqno £2.5.48,$$
 that $\varepsilon_{(\nu,\delta)} (m' +1) =\ell' +1$ and that 
 $(\mu_{(\nu,\delta)})_{m'+1}$ from $\frak t_{(\nu,\delta) \circ
(\mu,\varepsilon)} (\ell'\!+\!1)$ to $\frak r_{(\nu,\delta)}(m'\!+\!1)$  is determined by~$\mu_{m'+1}$ and
  $\frak t \big(\ell'\! +\!1 \bullet \varepsilon (m' \!+\!1)\big)\,;$ moreover, we consider the corresponding restriction
$$(\mu'_{(\nu,\delta)},\varepsilon'_{(\nu,\delta)}) : (\frak t'_{(\nu,\delta) \circ
(\mu,\varepsilon)},\Delta_{\ell' +1})\too (\frak r'_{(\nu,\delta)} ,\Delta_{m' +1})
\eqno £2.5.49\phantom{.}$$
which obviously fulfills
$$({\rm id}_{\frak r'_{(\nu,\delta)}},\iota^m_{m'})\circ
(\mu_{(\nu,\delta)},\varepsilon_{(\nu,\delta)}) =  (\mu'_{(\nu,\delta)},\varepsilon'_{(\nu,\delta)})
\circ ({\rm id}_{\frak t'_{(\nu,\delta)\circ (\mu,\varepsilon)}},\iota^\ell_{\ell'})
\eqno £2.5.50.$$

\smallskip
Now, it is easily checked that the composition $(\nu',\delta')\circ
(\mu'_{(\nu,\delta)},\varepsilon'_{(\nu,\delta)})$ coincides with the corresponding morphism~£2.5.43 for the $\ch^*(\F^{^\frak Y})\-$morphism $(\nu,\delta)\circ (\mu,\varepsilon)$ and therefore, by the very definition~£2.5.45, we have
$$\eqalign{\widehat\aut_{\F^{^\frak Y}}&\big((\nu,\delta)\circ (\mu,\varepsilon)\big)\cr
&=  \widehat\aut_{\F^{^\frak Y}} \big((\nu',\delta')\circ (\mu'_{(\nu,\delta)}, \varepsilon'_{(\nu,\delta)})\big) \circ
\widehat\aut_{\F^{^\frak Y}} \big({\rm id}_{\frak t'_{(\nu,\delta) \circ (\mu,\varepsilon)}},\iota^\ell_{\ell'}\big)\cr}
\eqno £2.5.51;$$
but, since $\vert R_{(\nu,\delta)}\vert/\vert \frak q (n)\vert < \vert\frak t (\ell)\vert/
\vert\frak q (n)\vert\,,$ it follows from the induction hypo-thesis that
$$\widehat\aut_{\F_{\!\rm nc}} \big((\nu',\delta') \circ
(\mu'_{(\nu,\delta)},\varepsilon'_{(\nu,\delta)})\big) = 
\widehat\aut_{\F^{^\frak Y}}(\nu',\delta')\circ\widehat\aut_{\F^{^\frak Y}}(\mu'_{(\nu,\delta)},\varepsilon'_{(\nu,\delta)})
\eqno £2.5.52;$$
similarly, since we have $\vert \frak t (\ell)\vert/\vert R_{(\nu,\delta)}\vert < 
\vert\frak t (\ell)\vert/\vert\frak q (n)\vert$ and
$$\widehat\aut_{\F^{^\frak Y}}(\mu_{(\nu,\delta)},\varepsilon_{(\nu,\delta)}) = \widehat\aut_{\F^{^\frak Y}} (\mu,\varepsilon)
\eqno £2.5.53,$$ 
we still get
$$\eqalign{&\widehat\aut_{\F^{^\frak Y}}\big((\nu,\delta)\circ (\mu,\varepsilon)\big)\cr 
& = \widehat\aut_{\F^{^\frak Y}}(\nu',\delta') \circ \widehat\aut_{\F^{^\frak Y}}
(\mu'_{(\nu,\delta)},\varepsilon'_{(\nu,\delta)}) \circ \widehat\aut_{\F^{^\frak Y}} ({\rm id}_{\frak t'_{(\nu,\delta)\circ (\mu,\varepsilon)}},\iota^\ell_{\ell'})\cr
&= \widehat\aut_{\F^{^\frak Y}}(\nu',\delta') \circ \widehat\aut_{\F^{^\frak Y}}
\big(({\rm id}_{\frak r'_{(\nu,\delta)}},\iota^m_{m'})\circ
(\mu_{(\nu,\delta)},\varepsilon_{(\nu,\delta)})\big)\cr  
&= \widehat\aut_{\F^{^\frak Y}}(\nu,\delta)\circ\widehat\aut_{\F^{^\frak Y}}
(\mu,\varepsilon)\,.\cr}
\eqno £2.5.54.$$

\smallskip
 Finally, we may assume that the image of $\frak r\big(\delta (n)\big)$  by 
 $\frak r\big(\delta (n)\bullet m\big)$ is normal in~$\frak r (m)\,,$ so that the image of $\frak t\big((\varepsilon\circ\delta)(n)\big)$ by
$\frak t\big((\varepsilon\circ\delta)(n)\!\bullet\!\varepsilon(m)\big)$ is normal 
in~$\frak t\big(\varepsilon(m)\big)\,;$ in particular, denoting by $\ell'$ the maximal element in~$\Delta_\ell - \Delta_{(\varepsilon\circ\delta) (n) -1}$ such that  the image of $\frak t\big((\varepsilon\circ\delta) (n)\big)$ 
by  $\frak t \big((\varepsilon\circ\delta) (n)\bullet \ell'\big)$ is normal 
in~$\frak t (\ell')\,,$ we have $\varepsilon (m)\le \ell'\,.$ If $\ell' = \ell$ then, by £2.5.41, we may assume that the image of $\frak t\big(\varepsilon (m)\big)$ is not normal 
in~$\frak t (\ell)$  and, denoting by~$\ell''\ge \varepsilon (m)$ the maximal element
 in~$\Delta_\ell$ such that  the image of 
$\frak t\big(\varepsilon (m)\big)$ by $\frak t (\varepsilon (m)\!\bullet\! \ell'')$ is normal
in~$\frak r (\ell'')\,,$ by our very definition (cf.~£2.5.45) we have
$$\widehat\aut_{\F^{^\frak Y}}(\mu,\varepsilon) =
\widehat\aut_{\F^{^\frak Y}} (\mu',\varepsilon') \circ\widehat\aut_{\F^{^\frak Y}}
({\rm id}_{\frak t'_{(\mu,\varepsilon)}},\iota^\ell_{\ell''})
\eqno £2.5.55;$$
but, according to equality~£2.5.41, we have
$$\widehat\aut_{\F^{^\frak Y}}(\nu,\delta)\circ \widehat\aut_{\F^{^\frak Y}} (\mu',\varepsilon') = \widehat\aut_{\F^{^\frak Y}}\big((\nu,\delta)\circ (\mu',\varepsilon')\big)
\eqno £2.5.56;$$
hence, since in the compositions of $(\nu,\delta)$ with $(\mu,\varepsilon)$ and
of $\big((\nu,\delta)\circ (\mu',\varepsilon')\big)$ with $ ({\rm id}_{\frak t'_{(\mu,\varepsilon)}},\iota^\ell_{\ell''})$ the first induction indices coincide with each other and the second ones strictly decreasse, it follows from the induction hypothesis that
$$\eqalign{\widehat\aut_{\F^{^\frak Y}}(\nu,\delta)&\circ
 \widehat\aut_{\F^{^\frak Y}} (\mu,\varepsilon)\cr
&= \widehat\aut_{\F^{^\frak Y}}(\nu,\delta)\circ \widehat\aut_{\F^{^\frak Y}}(\mu',\varepsilon')
\circ\widehat\aut_{\F^{^\frak Y}} ({\rm id}_{\frak t'_{(\mu,\varepsilon)}},\iota^\ell_{\ell''})
\cr
&= \widehat\aut_{\F^{^\frak Y}}\big((\nu,\delta)\circ (\mu',\varepsilon')\big)
\circ\widehat\aut_{\F^{^\frak Y}} ({\rm id}_{\frak
t'_{(\mu,\varepsilon)}},\iota^\ell_{\ell''})\cr
&= \widehat\aut_{\F^{^\frak Y}}\big((\nu,\delta)\circ (\mu,\varepsilon)\big)\cr}
\eqno £2.5.57.$$

\smallskip
In any case, we have a $\ch^*(\F^{^\frak Y})\-$morphism
$$(\mu'_{(\nu,\delta)},\varepsilon'_{(\nu,\delta)}) : (\frak t'_{(\nu,\delta) \circ
(\mu,\varepsilon)},\Delta_{\ell' +1})\too (\frak r ,\Delta_m)
\eqno £2.5.58\phantom{.}$$
fulfilling
$$(\mu'_{(\nu,\delta)},\varepsilon'_{(\nu,\delta)})\circ 
({\rm id}_{\frak t'_{(\nu,\delta) \circ (\mu,\varepsilon)}},\iota^\ell_{\ell'}) =
(\mu,\varepsilon) \circ ({\rm id}_\frak t, \delta^\ell_{\ell'+1})
\eqno £2.5.59;$$
 as above, it is easily checked that the composition $(\nu,\delta)\circ
(\mu'_{(\nu,\delta)},\varepsilon'_{(\nu,\delta)})$ coincides with the corresponding morphism~£2.5.43 for the $\ch^*(\F^{^\frak Y})\-$morphism $(\nu,\delta)\circ (\mu,\varepsilon)$ and therefore, by the very definition~£2.5.45, we have
$$\eqalign{\widehat\aut_{\F^{^\frak Y}}&\big((\nu,\delta)\circ (\mu,\varepsilon)\big)\cr
&= \widehat\aut_{\F^{^\frak Y}} \big((\nu,\delta)\circ (\mu'_{(\nu,\delta)},\varepsilon'_{(\nu,\delta)})\big) \circ
\widehat\aut_{\F^{^\frak Y}} \big({\rm id}_{\frak t'_{(\nu,\delta) \circ (\mu,\varepsilon)}},\iota^\ell_{\ell'}\big)\cr}
\eqno £2.5.60;$$
since  $\widehat\aut_{\F^{^{\rm nc}}}({\rm id}_\frak t,\delta^\ell_{\ell'+1}) = {\rm id}_{\hat\F (\frak t)}$ and we may assume that $\ell'\not= \ell$ ,
it follows from the induction hypothesis applied to the composition of $(\nu,\delta)$
with $(\mu'_{(\nu,\delta)},\varepsilon'_{(\nu,\delta)})$ that
$$\widehat\aut_{\F^{^\frak Y}} \big((\nu,\delta)\circ
(\mu'_{(\nu,\delta)},\varepsilon'_{(\nu,\delta)})\big) = \widehat\aut_{\F^{^\frak Y}}
(\nu,\delta)\circ \widehat\aut_{\F^{^\frak Y}}(\mu'_{(\nu,\delta)},\varepsilon'_{(\nu,\delta)})
\eqno £2.5.61;$$
moreover, if  $\vert \frak q (n)\vert < \vert \frak r (m)\vert\,,$ we can apply the induction hypothesis to both members of equality~£2.5.59 and then  we get
$$\widehat\aut_{\F^{^\frak Y}}(\mu'_{(\nu,\delta)},\varepsilon'_{(\nu,\delta)})\circ 
\widehat\aut_{\F^{^\frak Y}} \big({\rm id}_{\frak t'_{(\nu,\delta) \circ
(\mu,\varepsilon)}},\iota^\ell_{\ell'}\big)
 = \widehat\aut_{\F^{^\frak Y}}(\mu,\varepsilon)
 \eqno £2.5.62.$$
Consequently, once again we have
$$\widehat\aut_{\F_{\!\rm nc}}\big((\nu,\delta)\circ (\mu,\varepsilon)\big) = 
\widehat\aut_{\F_{\!\rm nc}} (\nu,\delta)\circ \widehat\aut_{\F_{\!\rm nc}}(\mu,\varepsilon)
 \eqno £2.5.63.$$

 \smallskip
 If $\vert \frak q (n)\vert = \vert \frak r (m)\vert$ then it follows from the  definitions of 
 $\widehat\aut_{\F^{^\frak Y}}(\mu,\varepsilon)$ and of $\widehat\aut_{\F^{^\frak Y}}\big((\nu,\delta)\circ (\mu,\varepsilon)\big)$ (cf.~£2.5.45) that $\ell'$ coincides with both induction  indices, that we get $\frak t'_{(\mu,\varepsilon)} = \frak t'_{(\nu,\delta)\circ  (\mu,\varepsilon)}$ and that the homomorphism~£2.5.43
 $$(\frak t'_{(\nu,\delta)\circ  (\mu,\varepsilon)},\Delta_{\ell' +1})
 \too (\frak q,\Delta_n)
 \eqno £2.5.64\phantom{.}$$
 corresponding to the  composition $(\nu,\delta)\circ (\mu,\varepsilon)$ coincides
 with $(\nu,\delta)\circ (\mu',\varepsilon')\,;$ at this point, we can apply 
 equality~£2.5.41 to obtain
 $$\widehat\aut_{\F^{^\frak Y}}(\nu,\delta)\circ\widehat\aut_{\F^{^\frak Y}}
 (\mu',\varepsilon') =\widehat\aut_{\F^{^\frak Y}}\big((\nu,\delta)\circ (\mu',\varepsilon')\big)
\eqno £2.5.65;$$
then, composing this equality with $\widehat\aut_{\F^{^\frak Y}}
({\rm id}_{\frak t'_{(\mu,\varepsilon)}},\iota^\ell_{\ell'})\,,$ from definition 
£2.5.45 we get
$$\widehat\aut_{\F^{^\frak Y}}(\nu,\delta)\circ\widehat\aut_{\F^{^\frak Y}}
 (\mu,\varepsilon) =\widehat\aut_{\F^{^\frak Y}}\big((\nu,\delta)\circ (\mu,\varepsilon)\big)
\eqno £2.5.66.$$
We are done.

\bigskip
\noindent
{\bf Theorem £2.6.[6, Theorem~£2.5]}\phantom{.} {\it Any functor $\widehat\aut_{\F^{^{\rm rd}}}$
lifting $\aut_{\F^{^{\rm rd}}}$ to the category $k^*\-\Gr$ can be extended to a unique folder structure of $\F\,.$\/}

\bigskip
\noindent
{\bf Theorem~£2.7.[5, Theorem~11.32]} {\it  The {\it Frobenius $P\-$category\/} $\F_{\!(b,G)}$ associated with~a block $b$ of a finite group $G$ has a
unique isomorphism class of folded structures  admitting
 a {\it $k^*\-$group isomorphism\/}
$$\widehat\aut_{\F_{\!(b,G)}^{^{\rm sc}}} (\frak q_Q)  \cong
\hat N_G(Q,f)/C_G (Q)
\eqno £2.7.1$$
for any $\F_{\!(b,G)}\-$selfcentralizing subgroup $Q$ of $P\,.$\/}

\medskip
£2.8. An obvious way for getting  a {\it folded structure\/} of $\F$
 is to start with a {\it regular central 
$k^*\-$extension\/} $\hat\F^{^{\rm sc}}$ of~$\F^{^{\rm sc}}\,;$ indeed, in this case  it follows again from [5, Proposition~A2.10] that we have a canonical functor
$$\aut_{\hat\F^{^{\rm sc}}} : \ch^*(\hat\F^{^{\rm sc}})\too k^*\-\Gr
\eqno £2.8.1$$
\noindent
mapping  any $\hat\F^{^{\rm sc}}\!\-$chain 
$\hat\frak q\,\colon \Delta_n\to \hat\F^{^{\rm sc}}$ to the stabilizer 
$\hat\F^{^{\rm sc}}\! (\frak q)$ in $\hat\F^{^{\rm sc}}\!
\big(\frak q (n)\big)$ of all the subgroups ${\rm Im}\big(\frak q(i\bullet n)\big)
$ for  $i\in \Delta_n\,,$ where $\frak q\,\colon \Delta_n\to \F^{^{\rm sc}}$
denotes the corresponding $\F^{^{\rm sc}}\!\-$chain; then, this functor
factorizes throughout a {\it folder structure\/} of $\F$
$$\widehat\aut_{\F^{^{\rm sc}}} : \ch^*(\F^{^{\rm sc}})\too k^*\-\Gr
\eqno £2.8.2.$$
Conversely, our main purpose here is to prove  that any {\it folder structure\/} of $\F$ comes from  a {\it regular central $k^*\-$extension\/} $\hat\F^{^{\rm sc}}$ 
of~$\F^{^{\rm sc}}$; consequently, once this result was obtained, to
 consider a {\it folded Frobenius $P\-$category\/} is equivalent to consider a pair 
 $(\F,\hat\F^{^{\rm sc}})$ formed by a {\it Frobenius  $P\-$category\/} $\F$ and by a {\it regular central $k^*\-$extension\/} 
$\hat\F^{^{\rm sc}}$ of $\F^{^{\rm sc}}\,.$ 
\eject

\medskip
£2.9. On the other hand,  in [1], [2] , [7] and [8]  it has been recently proved that there exists a unique {\it perfect $\F^{^{\rm sc}}\!\-$locality\/} $\P^{^{\rm sc}}$ [5,~17.4 and~17.13]. More explicitly, denote by $\T^{^{\rm sc}}_P$ the category  where the objects are all the {\it $\F\-$self-centralizing\/} subgroups of $P$ and,  for a pair of $\F\-$selfcentralizing subgroups $Q$ and $R$ of~$P\,,$ the set of morphisms from $R$ to $Q$ is the {\it $P\-$transporter\/} $T_P (R,Q)\,,$  the composition being induced by the product in $P\,;$ then [8,~\S4]
\smallskip
\noindent
£2.9.1\quad {\it there is a unique Abelian extension 
$\pi^{_{\rm sc}}\,\colon\P^{^{\rm sc}}\to \F^{^{\rm sc}}$ of  
$\F^{^{\rm sc}}$ endowed with a functor $\tau^{_{\rm sc}}\,\colon 
\T^{^{\rm sc}}_P \to \P^{^{\rm sc}}$  in such a way that the composition 
$\pi^{_{\rm sc}}\circ \tau^{_{\rm sc}}$ is the canonical functor defined by the conjugation in~$P\,,$ that $\P^{^{\rm sc}}\! (Q)$ is an {\it $\F\-$localizer\/} 
of $Q$ {\rm [5,~Theorem~18.6]} and that $Z(R)$ acts {\it regularly\/} over the {\it fibers\/} of the map $\P^{^{\rm sc}}(Q,R)\to \F^{^{\rm sc}}(Q,R)$
induced by $\pi^{_{\rm sc}}$ {\rm [5,~17.7]}, for any pair of $\F\-$selfcentralizing subgroups $Q$ and~$R$ of~$P\,.$\/}

\medskip
£2.10. Presently, the so-called 
{\it $\F\-$localizing  functor\/} considered in [6, 3.2.1] 
$$\loc_{\F^{^{\rm sc}}} : \ch^*(\F^{^{\rm sc}})\too \widetilde{\Loc}
\eqno £2.10.1\phantom{.}$$
is   just a {\it quotient\/} of the canonical functor [5, Proposition~A2.10]
$$\aut_{\P^{^{\rm sc}}} : \ch^* (\P^{^{\rm sc}})\too \Gr
\eqno £2.10.2.$$
Moreover, any {\it regular central $k^*\-$extension \/} $\hat\F^{^{\rm nc}}$ 
of $\F^{^{\rm nc}}$  determines {\it via\/}~$\pi^{_{\rm nc}}$ a {\it regular central 
$k^*\-$extension \/} $\hat\P^{^{\rm nc}}$ of $\P^{^{\rm nc}}\,;$ 
then,  the corresponding functor
$$\widehat{\loc}_{\F^{^{\rm sc}}} : \ch^* (\F^{^{\rm sc}})\too k^*\-\widetilde\Loc
\eqno £2.10.3\phantom{.}$$
considered in [6, 3.3.1] is   just a {\it quotient\/} of the obvious canonical functor  [5, Proposition~A2.10]
$$\aut_{\hat\P^{^{\rm sc}}} : \ch^* (\hat\P^{^{\rm sc}})\too k^*\-\Gr
\eqno £2.10.4.$$
Actually, it is clear that $\pi^{_{\rm nc}}$ induces an {\it equivalence\/} between
the so-called  {\it exterior quotients\/} $\tilde\F^{^{\rm nc}}$ of 
$\F^{^{\rm nc}}$ and $\tilde\P^{^{\rm nc}}$ of $\P^{^{\rm sc}}$  [5,~1.3];
that is to say, the quotients of $\F^{^{\rm sc}}$ and $\P^{^{\rm nc}}$ by the 
{\it inner automorphisms\/} of the objects are just isomorphic and, in particular, the 
{\it regular central $k^*\-$extensions\/} of 
$\tilde\F^{^{\rm nc}}\,,$ $\F^{^{\rm nc}}$ and $\P^{^{\rm nc}}$ are
clearly in bijective correspondence. In particular, a {\it folder structure\/} in $\F$ is equivalent to a functor
$$\widehat\aut_{\P^{^{\rm nc}}} : \ch^*(\P^{^{\rm nc}})\too k^*\-\Gr
\eqno £2.10.5\phantom{.}$$
lifting the canonical functor  $\aut_{\P^{^{\rm nc}}}\,.$

 \bigskip
\noindent
{\bf £3. Regular central $k^*\-$extensions of $\F^{^{\rm sc}}$}

\medskip
£3.1. Let $(\F,\widehat\aut_{\F^{^{\rm sc}}})$ be a {\it folded Frobenius 
$P\-$category\/} (cf.~£2.4) and denote by $\P$ and $\P^{^{\rm sc}}$ the
respective {\it perfect\/} $\F\-$ and $\F^{^{\rm sc}}\-${\it localities\/} [7,~\S6 and~\S7] and by 
$\pi\,\colon \P\to \F$ and $\tau\,\colon \T_P\to \P$ the {\it structural functors} [5,~17.3]. Our main prupose is to show that $(\F,\widehat\aut_{\F^{^{\rm sc}}})$ or, equivalently, 
$(\P,\widehat\aut_{\P^{^{\rm sc}}})$ (cf. £2.10.5) is determined by a
{\it regular central $k^*\-$extension\/}~$\hat\P^{^{\rm sc}}$ of 
$\P^{^{\rm sc}}$; we choose to work on $\P^{^{\rm sc}}$ rather than on 
$\F^{^{\rm sc}}$, which is equivalent as mentioned above, since in $\P^{^{\rm sc}}$all the morphisms are monomorphisms and epimorphisms [5, Proposition~24.2].
\eject

\medskip
£3.2. In particular, if $Q$ and $Q'$ are $\F\-$isomorphic $\F\-$selfcentralizing subgroups  of~$P\,,$ for any pair of  $\F\-$selfcentralizing  subgroups $R$ of $Q$ and $R'$ of~$Q'$ condition~£2.1.1 in $\F$ induces an injective {\it restriction\/} map
$$r^{Q',Q}_{R',R} : \P (Q',Q)_{R',R}\too \P (R',R)
\eqno £3.2.1\phantom{.}$$ 
where $\P (Q',Q)_{R',R}$ denotes the set of $x\in \P (Q',Q)$ such that 
$\pi_{Q',Q} (x)$ maps $R$ on $R'$; in particular, we may identify the stabilizer $\P (Q)_R$ of $R$ in $\P (Q)$ with a subgroup of $\P (R)\,.$ First of all, note the following consequence of condition~£2.1.3.

\bigskip
\noindent
{\bf Lemma~£3.3.} {\it With the notation above, assume that $R$ and $R'$
are $\F\-$isomor-phic and fully normalized in $\F\,;$  set $N = N_P(R)$ and $N' = N_P (R')\,.$ Then the restriction map and the composition induce a bijection
$$\P (N',N)_{R',R}\times_{\P (N)_R} \P (R) \cong \P (R',R)
\eqno £3.3.1.$$\/}

\par
\noindent
{\bf Proof:} It is clear that, for any $x\in \P (N',N)_{R',R}$ and any $s\in \P (R)$,~the composition $r^{N',N}_{R',R}(x)\. s$ belongs to $\P (R',R)\,;$ moreover,
 for any $y\in \P (N',N)_{R',R}$ and any  $t\in \P (R)$ such that $r^{N',N}_{R',R}(y)\. t =
 r^{N',N}_{R',R}(x)\. s\,,$  we clearly have that  $r^{N,N}_{R,R}(x^{-1}\.y) = s\.t^{-1}$
 which implies that $x^{-1}\.y$ belongs to $\P (N)_R\,;$ consequently, the pairs
 $(x,s)$ and $(y,t)$ have the same image in the quotient set
 $$\P (N',N)_{R',R}\times_{\P (N)_R} \P (R) = 
 \big(\P (N',N)_{R',R}\times \P (R)\big)\big/\P (N)_R
 \eqno £3.3.2.$$
 
 \smallskip
 Conversely, any $x\in \P (R',R)$ induces by conjugation a group isomorphism
 $\P(R)\cong \P(R')\,;$ then, since $\tau_R(N)$ and $\tau_{R'}(N')$ are respective Sylow  $p\-$subgroups of $\P (N)$ and $\P (N')$ [5,~2.11.4], there is $s\in \P(R)$
 such that the isomorphism $\P(R)\cong \P(R')$ induced by $x\.s$ sends
 $\tau_R(N)$ onto $\tau_{R'}(N')\,;$ at this point, it follows from condition~£2.1.3
 that there is $y\in \P (N',N)$ such that $r^{N',N}_{R',R}(y) = x\. s\,,$
 so that $y$ belongs to $\P (N',N)_{R',R}$ and $x$ is the image of the pair 
 $(y,s^{-1})\,.$

 \medskip
 £3.4. In order to discuss the uniqueness of the announced 
 {\it $k^*\-$category\/}~$\hat\P^{^{\rm sc}}$, note that the 
 {\it coherent $\F^{^{\rm sc}}\!\-$locality structure} of $\P^{^{\rm sc}}$ [5,~17.9] can be lifted to a  {\it coherent $\F^{^{\rm sc}}\!\-$locality structure} of $\hat\P^{^{\rm sc}}$. 
More precisely, let us consider  a nonempty set $\frak X$ of  $\F\-$selfcentralizing  subgroups of $P$ which contains any subgroup of $P$ admitting an $\F\-$morphism from some subgroup in~$\frak X\,,$ and respectively denote by 
$\T^{^{\frak X}}_P\,,$ $\F^{^\frak X}$ and $\P^{^\frak X}$ the {\it full\/} subcategories of $\T^{^{\rm sc}}_P\,,$  $\F^{^{\rm sc}}$ and 
$\P^{^{\rm sc}}$ over $\frak X$ as the  set of objects; we actually will prove that there exists an essentially unique regular central  $k^*\-$extension $\hat\P^{^\frak X}$ of $\P^{^\frak X}$
inducing  the obvious restricted functor (cf.~£3.1)
$$\widehat\aut_{\P^{^\frak X}} : \ch^*(\P^{^\frak X})\too k^*\-\Gr
\eqno £3.4.1;$$
 first of all, we claim that the  {\it coherent $\F^{^{\frak X}}\!\-$locality structure} of  $\P^{^{\frak X}}$ [5,~17.9] can be lifted to a  
{\it coherent $\F^{^{\frak X}}\!\-$locality structure} of $\hat\P^{^{\frak X}}$.

\bigskip
\noindent
{\bf Proposition~£3.5.} {\it With the notation above, the first structural functor 
 $\tau^{_{\frak X}}\,\colon \T^{^{\frak X}}_P \to \P^{^{\frak X}}$ can be lifted to
 a functor $\hat\tau^{_{\frak X}}\,\colon \T^{^{\frak X}}_P \to 
 \hat\P^{^{\frak X}}$ and such a lifting fulfills 
 $$\hat x\.\hat\tau^{_{\frak X}}_R (v) = 
 \hat\tau^{_{\frak X}}_Q\Big(\big(\pi_{Q,R} (x )\big)(v)\Big)\.\hat x 
\eqno £3.5.1\phantom{.}$$ 
for any pair of subgroups $Q$ and $R$ in $\frak X\,,$ any $x\in\P (Q,R)\,,$
any $\hat x\in \hat\P^{^{\frak X}}\! (Q,R)$ lifting $x$ and any $v\in R\,.$\/}

\medskip
\noindent
{\bf Proof:} We already know that $\tau_P\,\colon P\to \P (P)$
 is injective and thus, it can be uniquely lifted to an injective group homomrophism
 $\hat\tau_P\,\colon P\to \hat\P^{^{\frak X}}\! (P)\,;$ then,  choosing 
 $\hat\tau^{_{\frak X}}_{P,Q}(1)$ lifting $\tau_{P,Q}(1)$ in 
 $\hat\P^{^{\frak X}} (P,Q)$ for~any  subgroup $Q\not= P$  in~$\X\,,$ the functor $\hat\tau^{_{\frak X}}$ maps any  $\T^{^{\frak X}}_P\-$morphism 
 $u\,\colon R\to Q$ on the unique element 
 $\hat\tau^{_{\frak X}}_{Q,R}(u)$ in $\hat\P^{^{\frak X}}\! (Q,R)$ fulfilling
$$\hat\tau^{_{\frak X}}_{P,Q}(1)\.\hat\tau^{_{\frak X}}_{Q,R}(u) = 
\hat\tau_P (u)\.\hat\tau^{_{\frak X}}_{P,R}(1)
\eqno £3.5.2\phantom{.}$$
which makes sense since $u$ belongs to the {\it transporter\/} $T_P (R,Q)\,.$

\smallskip
With such a choice, $\hat\P^{^{\frak X}}\!$ becomes a {\it divisible $\F^{^{\frak X}}\!\-$locality\/} [5,~17.7], the {\it divi-sibility\/} being an easy consequence of the {\it divisibility\/} of $\P$ and of the {\it regularity\/} of the $k^*\-$extension
$\hat\P^{^{\frak X}}$; thus, our argument in [5,~Proposition~17.10] applies to
$\hat\P^{^{\frak X}}\!$ and therefore it suffices to prove condition~[5,~17.10.1];
but, note that for any $\hat x\in \hat\P^{^{\frak X}}\! (Q)$ the homomorphisms sending $v\in Q$ to  $\hat x\.\hat\tau^{_{\frak X}}_Q(v)\.\hat x^{-1}$ and 
to~$\hat\tau^{_{\frak X}}_Q\Big(\big(\pi_Q (x)\big)(v)\Big)$ lift the same 
group homomorphism from $Q$ to $\P(Q)$ and therefore they coincide with each other.

\medskip
£3.6. Note that, since a regular central  $k^*\-$extension $\hat\P^{^{\frak X}}$ of $\P^{^\frak X}$ endowed with a functor $\hat\tau^{_{\frak X}}\,\colon \T^{^{\frak X}}_P \to \hat\P^{^{\frak X}}$ lifting  the first structural functor  
$\tau^{_{\frak X}}\,\colon \T^{^{\frak X}}_P \to \P^{^{\frak X}}$ and fulfilling condition~£3.5.1 is actually a {\it coherent $\F^{^{\frak X}}\!\-$locality\/} 
[5,~£17.7], with the notation in~£3.2 above we also have an injective 
{\it $k^*\-$restriction\/} map
$$\hat r^{Q',Q}_{R',R} : \hat\P^{^{\frak X}}\! (Q',Q)_{R',R}\too 
\hat\P^{^{\frak X}}\! (R',R)
\eqno £3.6.1\phantom{.}$$ 
where $\hat\P^{^{\frak X}}\! (Q',Q)_{R',R}$ is the converse image of
$\P (Q',Q)_{R',R}$ in $\hat\P^{^{\frak X}}\! (Q',Q)\,.$

 \bigskip
 \noindent
 {\bf Theorem~£3.7.} {\it With the notation above, there exists a regular central  
 $k^*\-$ex-tension $\hat\P^{^{\rm sc}}$ of $\P^{^{\rm sc}}\,,$ unique up to 
 $k^*\-$equivalences, inducing the  folded Frobenius $P\-$category 
 $(\F,\widehat\aut_{\F^{^{\rm sc}}})\,.$\/}

\medskip
\noindent
{\bf Proof:} We choose a set $\frak X$ as above and, arguing by induction
on $\vert\frak X\vert\,,$ we will prove that there exists a regular central  
$k^*\-$extension $\hat\P^{^\frak X}$ of $\P^{^\frak X}$
inducing  the obvious restricted functor (cf.~£3.1)
$$\widehat\aut_{\F^{^\frak X}} : \ch^*(\F^{^\frak X})\too k^*\-\Gr
\eqno £3.7.1\phantom{.}$$
and that such a $\hat\P^{^\frak X}$ endowed with a lifting $\hat\tau^{_{\frak X}}\,\colon \T^{^{\frak X}}_P \to \hat\P^{^{\frak X}}$ of $\tau^{_{\frak X}}\,,$
which fulfills condition~£3.5.1, is unique up to $k^*\-$equivalences.

\smallskip
If $\frak X = \{P\}$ then $\P^{^\frak X}$ has just one object $P$ and its automorphism group is $\P (P)\,;$ then, the {\it folder structure\/} maps the
trivial {\it $\F^{^{\rm sc}}\-$chain\/} $\Delta_0\to \F^{^{\rm sc}}$
sending $0$ to $P$ on a $k^*\-$group $\hat\F (P)$ which, by restriction, determines a $k^*\-$group~$\hat\P (P)\,;$ that is to say, we get a 
$k^*\-$category $\hat\P^{^\frak X}$ with one object~$P$ and with the 
$k^*\-$group automorphism $\hat\P(P)\,,$ which clearly induces the 
corresponding functor~£3.7.1 again; the uniqueness is clear.

\smallskip
Otherwise,  choose a minimal element $U$ in $\frak X$ {\it fully normalized\/} in 
$\F$ and set 
$$\frak Y = \frak X - \{\theta(U)\mid \theta\in \F(P,U)\}
\eqno £3.7.2;$$
that is to say, according to our induction hypothesis,  there exists a regular central  $k^*\-$extension $\hat\P^{^\frak Y}$ of $\P^{^\frak Y}$
inducing  the obvious restricted functor (cf.~£3.1)
$$\widehat\aut_{\F^{^\frak Y}} : \ch^*(\F^{^\frak Y})\too k^*\-\Gr
\eqno £3.7.3.$$
and such a $k^*\-$category $\hat\P^{^\frak Y}$ endowed with a lifting 
$\hat\tau^{_{\frak Y}}\,\colon \T^{^{\frak Y}}_P \to \hat\P^{^{\frak Y}}$ of 
$\tau^{_{\frak Y}}$ which fulfills condition~£3.5.1 (cf.~Proposition~£3.5) is unique up to $k^*\-$isomorphisms.

\smallskip
If $N_\F (U) = \F$ [5,~Proposition~2.16], we also have $N_\P (U) = \P$
[5,~17.5] and then it is easily checked from~£3.2.1 that $\P^{^\frak X}$ actually
coincides with the category $\T^{^\frak X}_{\P (U)}$  where $\frak X$ is the set of objects and where,   for a pair of subgroups $Q$ and $R$ in~$\frak X\,,$ 
 the set of morphisms from $R$ to $Q$ is the {\it $\P (U)\-$transporter\/} 
 $$\T^{^\frak X}_{\P (U)} (Q,R) = \{x\in \P (U)\mid x\.\tau_U (R)\.x^{-1}\i
 \tau_U (Q)\}
 \eqno £3.7.4,$$
 the composition being defined by the product in $\P (U)\,;$
  but, once again,  the {\it folder structure\/} maps the
trivial {\it $\F^{^{\rm sc}}\-$chain\/} $\Delta_0\to \F^{^{\rm sc}}$
sending $0$ to $U$ on a $k^*\-$group $\hat\F (U)$ which, by restriction, determines a $k^*\-$group~$\hat\P (U)\,;$ hence, denoting by $\hat\tau_U (Q)$
and $\hat\tau_U (R)$ the finite $p\-$subgroups of~$\hat\P (U)$ respectively lifting
$\tau_U (Q)$ and $\tau_U (R)\,,$ we can consider the corresponding {\it transporter\/} in the $k^*\-$group 
$\hat\P (U)$ 
 $$\T^{^\frak X}_{\hat\P (U)} (Q,R) = \{\hat x\in \hat\P (U)\mid 
 \hat x\.\hat\tau_U (R)\. \hat x^{-1}\i  \hat\tau_U (Q)\}
 \eqno £3.7.5.$$

 \smallskip
 Now, it is clear that the $k^*\-$category $\T^{^\frak X}_{\hat\P (U)}$  where 
 $\frak X$ is the set of objects, where the obvious $k^*\-$set 
 $\T^{^\frak X}_{\hat\P (U)} (Q,R)$ is  the $k^*\-$set of morphisms from $R$ 
 to~$Q$ for any pair of subgroups $Q$ and $R$ in~$\frak X\,,$ and where the composition is defined~by the product in $\hat\P (U)$ determines a 
 {\it regular central  $k^*\-$extension\/} of~$\T^{^\frak X}_{\P (U)}  
 = \P^{^\frak X}$ together with an obvious lifting of $\tau^{_{\frak X}}\,,$
 which fulfills condition~£3.5.1.

 \smallskip
  On the other hand,  it is easily checked that  such a  {\it regular central 
 $k^*\-$ex-tension\/} $\hat\P^{^\frak X}\!$ is also {\it divisible\/} [5,~17.7] and
 therefore that,  for any pair of subgroups $Q$ and $R$ in~$\frak X\,,$
 as in~£3.2.1 above we get a {\it restriction\/} $k^*\-$set homomorphism
 $$\hat\P^{^\frak X} (Q\.U,R\.U)\too \hat\P^{^\frak X} (U)
 \eqno £3.7.6\phantom{.}$$
 which is always injective; moreover, since we have $N_\P (U) = \P\,,$
 always by the divisibility of $\hat\P^{^\frak X}\!$ we get a $k^*\-$set isomomorphism
 $$\hat\P^{^\frak X} (Q\.U,R\.U)_{Q,R}\cong \hat\P^{^\frak X} (Q,R)
 \eqno £3.7.7\,.$$
 From these remarks, it is easily checked the uniqueness of $\hat\P^{^\frak X}$
 and the fact that this {\it $k^*\-$category\/} 
determines the restricted functor $\widehat\aut_{\F^{^\frak X}}\,.$

\smallskip
Otherwise recall that, according to [6,~3.1], for any subgroup $Q$ of $P$ fully normalized in $\F\,,$  our {\it folded Frobenius $P\-$category\/} induces a  {\it folded Frobenius  $N_P(Q)\-$category\/} 
$\big(N_\F (Q),\widehat\aut_{N_\F (Q)^{^{\rm sc}}}\big)$ where 
$$\widehat\aut_{N_\F (Q)^{^{\rm sc}}} : \ch^*(N_\F (Q)^{^{\rm sc}})
\too k^*\-\Gr
\eqno £3.7.8\phantom{.}$$ 
is the unique functor lifting $\aut_{N_\F (Q)^{^{\rm sc}}}$ and extending the restriction of $\widehat\aut_{\F^{^{\rm sc}}}$ to $N_\F (Q)^{^{\rm rd}}$
(cf.~Theorem~£2.6 and [6, Lemma 2.5]).

\smallskip
 Thus, if we have $N_\F (U) \not= \F\,,$ arguing by induction on the size of $\F\,,$ for any $V\in \frak X - \frak Y$ fully normalized in $\F$ we may assume that there exists a {\it regular central  $k^*\-$extension\/} $\widehat{N_\P (V)}^{_{\rm sc}}$ of 
$N_\P (V)^{^{\rm sc}}$ determining  $\widehat\aut_{N_\F (V)^{^{\rm sc}}}\,,$~and that such a {\it $k^*\-$category\/} $\widehat{N_\P (V)}^{_{\rm sc}}$, endowed with a lifting 
$\hat\tau^{_{V,\rm sc}}: \T^{^{\rm sc}}_{N_P(V)}\to\widehat{N_\P (V)}^{_{\rm sc}}$ of the first {\it structural functor\/} of $N_\F (V)^{^{\rm sc}}\!$ which fulfills condition~£3.5.1 (cf.~Pro-position~£3.5), is unique up to  $k^*\-$isomorphisms. Actually, we are only in-terested in the {\it full\/} 
$k^*\-$subcategory of $\widehat{N_\P (V)}^{_{\rm sc}}$ over  the set 
$N_\frak X (V)$ of subgroups in $\frak X$   contained in $N_P (V)$ and may assume that the lifting 
$$\hat\tau^{_{V,N_\frak Y (V)}} :  \T^{^{N_\frak Y (V)}}_{N_P(V)} \too
\widehat{N_\P (V)}^{_{N_\frak Y (V)}}
\eqno £3.7.9\phantom{.}$$
coincides with the restriction of $\hat\tau^{_{\frak Y}}\,;$ then, it follows from Proposition~£3.5 that we can identify $\widehat{N_\P (V)}^{_{N_\frak Y (V)}}$
with the {\it full\/} $k^*\-$subcategory of $\hat\P^{^{\frak Y}}$ over the set 
$N_\frak Y (V)\,.$

\smallskip
Moreover, setting $N = N_P (V)$ and considering the 
{\it $N_\F (V)^{^{\rm sc}}\-$chains\/}  $\frak q_V\,\colon 
\Delta_0\to N_\F (V)^{^{\rm sc}}\,,$   $\frak q_N\,\colon \Delta_0\to 
N_\F (V)^{^{\rm sc}}$ (cf.~£2.2) and  $\frak n\,\colon \Delta_1\to 
N_\F (V)^{^{\rm sc}}$  which  map $0$ on $V\,,$ $1$ on $N$ and 
$0\bullet 1$ on the inclusion of $V$ in $N\,,$ noted $\iota_V^N\,,$ and the obvious 
$\ch^*(N_\F (V)^{^{\rm sc}})\-$morphisms (cf.~£2.2)
$$(\id_V,\delta_1^0) : (\frak n,\Delta_1)\to (\frak q_V,\Delta_0)\qq
(\id_N,\delta_0^0) : (\frak n,\Delta_1)\to (\frak q_N,\Delta_0)
\eqno £3.7.10,$$
the functors $\widehat\aut_{N_\F (V)^{^{\rm sc}}}$ and 
$\widehat\aut_{\F^{^{\rm sc}}}$ send $\frak n\,,$  $\frak q_V$ and $\frak q_N$ to the same respective $k^*\-$groups $\hat\F (N)_V\,,$  $\hat\F (V)$ and
$\hat\F (N)\,,$ and they send the $\ch^*(N_\F (Q)^{^{\rm sc}})\-$mor-phisms 
$(\id_V,\delta_1^0)$ and $(\id_N,\delta_0^0)$ to the same respective $k^*\-$group  homomorphisms
$$\hat\F (N)_V\too \hat\F (V)\qq \hat\F (N)_V\too \hat\F (N)
\eqno £3.7.11;$$
note that the images of $\hat\F (N)_V$ are respectively 
$N_{\hat\F (V)}\big(\F_{\!N }(V)\big)$ and the stabilizer $\hat\F (N)_V$ of $V$ in 
$\hat\F (N)\,.$

\smallskip
Since $N$ belongs to $\frak Y\,,$ the restriction of $\hat\F (N)$ from $\F (N)$ to
$\P(N)$ ne-cessarily coincides with $\hat\P^{^\frak Y}\!(N)$ and therefore 
 the restriction of $\hat\F (N)_V$ from $\F (N)_V$ to $\P(N)_V$ also coincides
 with  the stabilizer $\hat\P^{^\frak Y} (N)_V$ of $V$ 
 in~$\hat\P^{^\frak Y}\! (N)\,.$ Then,  for any $V'\in \frak X - \frak Y$ fully normalized in $\F\,,$ setting $N' = N_P (V')$ and denoting by 
 $\hat\P^{^\frak Y}\! (N',N)_{V',V}$ the converse image of $\P (N',N)_{V',V}$ in
 $\hat\P^{^\frak Y}\! (N',N)$ and by $\hat\P^{^\frak X}(V)$  the restriction of 
 $\hat\F (V)$ from $\F (V)$ to $\P(V)\,,$ it is clear that 
 $\hat\P^{^\frak Y} (N)_V$ acts on the $k^*\-$set  $\hat\P^{^\frak Y}\! (N',N)_{V',V}$ by right-hand composition
 in $\hat\P^{^\frak Y}$; moreover, the left-hand homomorphism in~£3.7.10 induces a $k^*\-$group {\it injective\/} homomorphism form $\hat\P^{^\frak Y} (N)_V$ to 
 $\hat\P^{^\frak X} (V)\,;$ thus,  we are able to define the $k^*\-$set
$$\hat\P^{^\frak X}\! (V',V) = \hat\P^{^\frak Y}\! (N',N)_{V',V}
\times_{\hat\P^{^\frak Y} (N)_V} \hat\P^{^\frak X} \!(V) 
\eqno £3.7.12\phantom{.}$$
and then, from isomorphism~£3.3.1, we get a canonical map
$$\hat\P^{^\frak X}\! (V',V)\too \P (V',V)
\eqno £3.7.13.$$

\smallskip
Note that, in the case where $V' = V\,,$ our notation is coherent. Moreover, for another $V''\in \frak X - \frak Y$ fully normalized in $\F\,,$ setting $N'' = N_P (V'')$
and considering $\hat\P^{^\frak Y}\! (N'',N)_{V'',V}\,,$ 
$\hat\P^{^\frak Y}\! (N'',N')_{V'',V'}$ and $\hat\P^{^\frak X} \!(V')$ as above,
we also have the $k^*\-$sets
$$\eqalign{\hat\P^{^\frak X}\! (V'',V) &= \hat\P^{^\frak Y}\! (N'',N)_{V'',V}
\times_{\hat\P^{^\frak Y} (N)_V} \hat\P^{^\frak X} \!(V) \cr
\hat\P^{^\frak X}\! (V'',V') &= \hat\P^{^\frak Y}\! (N'',N')_{V'',V'}
\times_{\hat\P^{^\frak Y} (N')_{V'}} \hat\P^{^\frak X} \!(V') \cr}
\eqno £3.7.14\phantom{.}$$
and we claim that the composition in $\hat\P^{^\frak Y}$ and in the corresponding $k^*\-$groups induces a {\it $k^*\-$composition\/}
$$c^{_\frak X}_{V'',V',V} : \hat\P^{^\frak X}\! (V'',V') \times 
\hat\P^{^\frak X}\! (V',V)\too \hat\P^{^\frak X}\! (V'',V)
\eqno £3.7.15\phantom{.}$$
lifting the  composition in $\P$ {\it via\/} the canonical maps~£3.7.13.
\eject

\smallskip
First of all, {\it mutatis mutandis\/} denote by $\frak q_{V'}\,,$  
$\frak q_{N'}$ and $\frak n'\,,$  the analogous 
{\it $N_\F (V')^{^{\rm sc}}\!\-$chains\/} and by 
$(\id_{V'},\delta^0_1)$ and $(\id_{N'},\delta^0_1)$ the analogous
$\ch^*(N_\F (V')^{^{\rm sc}})\-$ morphisms, as in £3.7.10 above; it is clear that any $\F\-$morphism
$\varphi\,\colon N\to N'$ fulfilling $\varphi (V) = V'$ determines {\it natural 
isomorphisms\/} $\frak q_V\cong \frak q_{V'}\,,$ $\frak q_N\cong \frak q_{N'}$
and $\frak n\cong \frak n'$  which induce commutative 
$\ch^*(\F^{^{\rm sc}})\-$diagrams (cf.~£3.7.10)
$$\matrix{(\frak n',\Delta_1)&\too &(\frak q_{V'},\Delta_0)\cr
\wr\Vert&&\wr\Vert\cr
(\frak n,\Delta_1)&\too &(\frak q_V,\Delta_0)\cr}\qq
\matrix{(\frak n',\Delta_1)&\too &(\frak q_{N'},\Delta_0)\cr
\wr\Vert&&\wr\Vert\cr
(\frak n,\Delta_1)&\too &(\frak q_N,\Delta_0)\cr}
\eqno £3.7.16\,;$$
at this point, the functor $\widehat\aut_{\F^{^{\rm sc}}}$ sends these commutative $\ch^*(\F^{^{\rm sc}})\-$diagrams to the commutative 
 diagrams of $k^*\-$groups
$$\matrix{\hat\F (N')_{V'}&\too &\hat\F (V')\cr
\wr\Vert&&\wr\Vert\cr
\hat\F (N)_V&\too &\hat\F (V)\cr}\qq 
\matrix{\hat\F (N')_{V'}&\too &\hat\F (N')\cr
\wr\Vert&&\hskip-15pt{\scriptstyle \hat\frak f_\varphi}\hskip5pt\wr\Vert\cr
\hat\F (N)_V&\too &\hat\F (N)\cr}
\eqno £3.7.17.$$

\smallskip
Consequently, for any $x\in \P (N',N)_{V',V}$ lifting $\varphi$ we get  the commutative  diagrams of $k^*\-$groups
$$\matrix{\hat\P^{^\frak Y}\! (N')_{V'}&\too &\hat\P^{^\frak X}\! (V')\cr
\wr\Vert&&\hskip-15pt{\scriptstyle \hat\frak h_x}\hskip5pt\wr\!\Vert\cr
\hat\P^{^\frak Y}\! (N)_V&\too &\hat\P^{^\frak X}\! (V)\cr}\qq 
\matrix{\hat\P^{^\frak Y} (N')_{V'}&\too &\hat\P^{^\frak Y} (N')\cr
\wr\Vert&&\hskip-15pt{\scriptstyle \hat\frak g_x}\hskip5pt\wr\!\Vert\cr
\hat\P^{^\frak Y} (N)_V&\too &\hat\P^{^\frak Y} (N)\cr}
\eqno £3.7.18\phantom{.}$$
and note that the $k^*\-$group isomorphism $\hat\frak g_x$ has to be induced by the composition in $\hat\P^{^\frak Y}$ (cf.~£3.7.3); that is to say, for any 
$\hat x\in \hat\P^{^\frak Y}\! (N',N)_{V',V}$ lifting~$x$ and any 
$\hat s\in \hat\P^{^\frak Y} (N)\,,$ we actually have 
$\hat\frak g_x (\hat s) = \hat x\.\hat s\.\hat x^{-1}\,.$

\smallskip
We are ready to define the $k^*\-$composition $c^{_\frak X}_{V'',V',V}$ in~£3.7.15; any ele-ment in
$\hat\P^{^\frak X}\! (V',V)$ is the class $\,\overline{\!(\hat x,\hat s)\!}\,$ of some pair $(\hat x,\hat s)$ where $\hat x$ and $\hat s$ respectively belong to 
$\hat\P^{^\frak Y}\! (N',N)_{V',V}$ and to $\hat\P^{^\frak X}\! (V)\,;$
similarly, if $\,\overline{\!(\hat x',\hat s')\!}\,$ is an element 
of~$\hat\P^{^\frak X}\! (V'',V')\,,$ it is clear that, in the $k^*\-$category
 $\hat\P^{^\frak Y}\,,$ the composition $\hat x'\.\hat x$ makes sense and belongs
to $\hat\P^{^\frak Y}\! (N'',N)_{V'',V}\,;$ moreover, denoting by $x$ the image of
$\hat x$ in $\P (N',N)\,,$ we have  the $k^*\-$group isomorphism  
$\hat\frak h_x$ from $\hat\P^{^\frak X} (V)$ to $\hat\P^{^\frak X} (V')$
and therefore $(\hat\frak h_x)^{-1}(\hat s')$ belongs to 
$\hat\P^{^\frak X} (V)\,;$ then, we set
$$c^{_\frak X}_{V'',V',V} \big(\,\overline{\!(\hat x',\hat s')\!}\,,\,\overline{\!(\hat x,\hat s)\!}\,\big) = \,\overline{\!\big(\hat x'\.\hat x,(\hat\frak h_x)^{-1}(\hat s')\.\hat s\big)\!}\,
\eqno £3.7.19;$$
the compatibility with the action of $k^*$ is clear.

\smallskip
This makes sense since, for any $\hat t\in \hat\P^{^\frak Y} (N)_{V}$ and any
$\hat t'\in \hat\P^{^\frak Y} (N')_{V'}\,,$ denoting by $t$ the image of $\hat t$
in $\P (N)$ we get (cf.~£3.7.18)
$$\eqalign{(\hat x'\.\hat t')\.(\hat x\.\hat t) 
&= \hat x'\.\hat x\.(\hat\frak g_x)^{-1} (\hat t')\.\hat t \cr
(\hat\frak h_{x\.t})^{-1}(\hat t'^{-1}\.\hat s')\.(\hat t^{-1}\.\hat s)
&= \big((\hat\frak h_t)^{-1}\circ (\hat\frak h_x)^{-1}\big)(\hat t'^{-1}\.\hat s')   \.\hat t^{-1}\.\hat s\cr
&= (\hat\frak h_t)^{-1}\big((\hat\frak g_x)^{-1}(\hat t'^{-1})\.
(\hat\frak h_x)^{-1}(\hat s')   \big)\.\hat t^{-1}\.\hat s\cr
&=  \hat t^{-1}\.(\hat\frak g_x)^{-1}(\hat t'^{-1})
\.(\hat\frak h_x)^{-1}(\hat s')\.\hat s  \cr
&=  \big((\hat\frak g_x)^{-1} (\hat t')\.\hat t\big)^{-1}\.(\hat\frak h_x)^{-1}(\hat s')\.\hat s\cr}
\eqno £3.7.20.$$
The $k^*\-$composition is associative since, for any $V'''\in \frak X - \frak Y$ fully normalized in $\F$ and any element $\,\overline{\!(\hat x'',\hat s'')\!}\,$
in $\hat\P^{^\frak X}\! (V''',V'')\,,$ denoting by $x'$ the image of $\hat x'$
in $\P (N'',N')$ we obtain
$$\eqalign{c^{_\frak X}_{V''',V'',V} 
&\Big(\,\overline{\!(\hat x'',\hat s'')\!}\,, c^{_\frak X}_{V'',V',V} 
\big(\,\overline{\!(\hat x',\hat s')\!}\,,\,\overline{\!(\hat x,\hat s)\!}\,\big)\Big)\cr
&= c^{_\frak X}_{V''',V'',V} \Big(\,\overline{\!(\hat x'',\hat s'')\!}
\,, \,\overline{\!\big(\hat x'\.\hat x,(\hat\frak h_x)^{-1}
(\hat s')\.\hat s\big)\!}\,\Big)\cr
&= \,\overline{\!\Big(\hat x''\.(\hat x'\.\hat x),
(\hat\frak h_{x'\. x})^{-1}(\hat s'')\.\big((\hat\frak h_x)^{-1}
(\hat s')\.\hat s\big)\Big)\!}\,\cr
&= \,\overline{\!\Big((\hat x''\.\hat x')\.\hat x,
(\hat\frak h_x)^{-1} \big((\hat\frak h_{x'})^{-1}(\hat s'')\.\hat s'\big)\.\hat s\Big)\!}\,\cr
&= c^{_\frak X}_{V''',V',V} \Big(c^{_\frak X}_{V'',V'',V'} 
\big(\,\overline{\!(\hat x'',\hat s'')\!}\,,\,\overline{\!(\hat x',\hat s')\!}\,\big)
,\,\overline{\!(\hat x,\hat s)\!}\,\Big)\cr}
\eqno £3.7.21.$$

\smallskip
According to our definition of $\hat\P^{^\frak X}\!  (V',V)$ in~£3.7.12, the unity element of $\hat\P^{^\frak X}\! (V)$ defines a canonical $k^*\-$set homomorphism
$$\hat r^{N',N}_{V',V} : \hat\P^{^\frak Y}\!  (N',N)_{V',V}\too 
\hat\P^{^\frak X}\!  (V',V)
\eqno £3.7.22\phantom{.}$$
lifting $r^{N',N}_{V',V}\,.$ More generally, let $Q$ and $Q'$ be a pair of 
subgroups of $P$ respectively contained in $N$ and $N'\,,$ and strictly containing $V$ and $V'\,;$ we define as follows an injective $k^*\-$set homomorphism
$$\hat r^{Q',Q}_{V',V} : \hat\P^{^\frak Y}\!  (Q',Q)_{V',V}\too 
\hat\P^{^\frak X}\!  (V',V)
\eqno £3.7.23\phantom{.}$$
lifting the restriction map (cf. £3.2.1)
$$r^{Q',Q}_{V',V} : \P (Q',Q)_{V',V}\too \P (V',V)
\eqno £3.7.24.$$
If $\hat x\in \hat\P^{^\frak Y}\!(Q',Q)_{V',V}$ and $x$ denotes its image in 
$\P (Q',Q)_{V',V}\,,$ it follows from Lemma~£3.3 that 
$r^{Q',Q}_{V',V} (x) = r^{N',N}_{V',V} (y)\.z$ for suitable $y\in \P (N',N)_{V',V}$
and $z\in \P (V)\,;$ thus, setting $Q'' = \big(\pi_{N,N'}(y^{-1})\big)(Q')\i N\,,$ 
we get
$$z = r^{Q'',Q}_{V,V} \big(r^{N,N'}_{Q'',Q'} (y^{-1})\.x\big)
\eqno £3.7.25\phantom{.}$$
and therefore, setting $s = r^{N,N'}_{Q'',Q'} (y^{-1})\.x\,,$ by injectivity of $r^{Q',Q}_{V',V}$  (cf.~£3.2) we still get $x = r^{N',N}_{Q',Q''} (y)\.s$. 
\eject

\smallskip
Hence, choosing a lifting $\hat y$ of $y$ in 
$\hat\P^{^\frak Y}\!  (N',N)_{V',V}\,,$ in the $k^*\-$category 
$\hat\P^{^\frak Y}\!$ we have the restriction $\hat r^{N',N}_{Q',Q''} (\hat y)$
 (cf.~£3.6) as an element of $\hat\P^{^\frak Y}\!  (Q',Q'')_{V',V}\,;$
then, there is a unique lifting $\hat s$ of $s$ in 
$\hat\P^{^\frak Y}\!  (Q'',Q)_{V,V}$ fulfilling 
$\hat x = \hat r^{N',N}_{Q',Q''} (\hat y)\.\hat s\,.$ Moreover, since 
$\widehat{N_\P (V)}^{_{N_\frak Y (V)}}$ can be identified with the {\it full\/} 
$k^*\-$subcategory of $\hat\P^{^{\frak Y}}$ over the set $N_\frak Y (V)\,,$ actually $\hat s$ can be identified with an element of
$\widehat{N_\P (V)}^{_{\rm sc}}\! (Q'',Q)$ stabilizing $V$ and therefore in the $k^*\-$category $\widehat{N_\P (V)}^{_{N_\frak X (V)}}$ we have the restriction
$\hat r_{V,V}^{Q'',Q}(\hat s)$  (cf.~£3.6) lifting $z$ to 
$\widehat{N_\P (V)}^{_{N_\frak X (V)}} (V)$ which coincides with
$\hat\P^{^\frak X}\!  (V)$ since we have
$$\widehat{N_\F (V)}^{_{\rm sc}}\! (V) = \widehat\aut_{N_\F (V)^{^{\rm sc}}} (\frak q_V) = \widehat\aut_{\F^{^{\rm sc}}} (\frak q_V) = \hat\F (V)
\eqno £3.7.26.$$
Then, we define (cf.~£3.7.12)
$$\hat r^{Q',Q}_{V',V}(\hat x) = 
\,\overline{\!\big(\hat y,\hat r_{V,V}^{Q'',Q}(\hat s)\,\big)\!}
\eqno £3.7.27;$$
it is independent of our choice of $y\in \P (N',N)_{V',V}$  since, for another decompostion $r^{Q',Q}_{V',V} (x) = r^{N',N}_{V',V} (y')\.z'\,,$ we actually have
$y' = y\.t$ and $z' = r^{N}_{V} (t^{-1})\.z$ for some $t\in \P (N)_V\,;$
thus, setting  $Q''' = \big(\pi_{N}(t^{-1})\big)(Q'')\,,$ once again an element 
$\hat t$ of $\hat\P^{^\frak Y}\! (N)_V$ lifting~$t$  can be identified  with an element of $\widehat{N_\P (V)}^{_{\rm sc}}\! (N)$ stabilizing $V$  and we also obtain 
$$\hat x = \hat r^{N',N}_{Q',Q''} (\hat y)\.\hat s = \big( \hat r^{N',N}_{Q',Q'''} (\hat y\.\hat t)\big)\.\big(\hat r^{N,N}_{Q''',Q''}(\hat t^{-1})\.\hat s\big)
\eqno £3.7.28;$$
but, the pairs $\big(\hat y,\hat r_{V,V}^{Q'',Q}(\hat s)\,\big)$ and 
$\big(\hat y\.\hat t,\hat r_{V,V}^{Q''',Q}(\hat r^{N,N}_{Q''',Q''}
(\hat t^{-1})\.\hat s)\,\big)$ have the same class in 
$\hat\P^{^\frak X}\!  (V',V)\,.$

\smallskip
At present, if $R$ and $R'$ are a pair of subgroups of $P$ respectively contained in $Q$ and $Q'\,,$ and strictly containing $V$ and $V'\,,$ we claim that the corresponding restriction $\hat r^{R',R}_{V',V}$ agree with 
$\hat r^{Q',Q}_{V',V}\,;$ if $\hat x\in \hat\P^{^\frak Y}\!  (Q',Q)_{V',V}$ has
an image in $\F (Q',Q)$ mapping $R$ on $R'\,,$ it follows from £3.6 above that we have the restriction  $\hat r_{R',R}^{Q',Q}(\hat x)$ in 
$\hat\P^{^\frak Y}\!  (R',R)_{V',V}$ and we claim that
$$\hat r^{R',R}_{V',V} \big(\hat r_{R',R}^{Q',Q}(\hat x)\big) = \hat r^{Q',Q}_{V',V}(\hat x)
\eqno £3.7.29;$$
indeed, with the notation above we may assume that 
$\hat x = \hat r^{N',N}_{Q',Q''} (\hat y)\.\hat s\,;$ then, setting 
$R'' = \big(\pi_{N,N'}(y^{-1})\big)(R')\i N\,,$  we clearly have
$$\hat r_{R',R}^{Q',Q}(\hat x) = \hat r^{N',N}_{R',R''} (\hat y)\.
\hat r_{R'',R}^{Q'',Q}(\hat s)
\eqno £3.7.30;$$
consequently, considering the set $N_\frak X (V)$ defined above,  since the restriction  in the $k^*\-$category 
$\widehat{N_\P (V)}^{_{N_\frak X (V)}}$ is transitive (cf.~£3.6), we clearly obtain
$$\hat r^{R',R}_{V',V} \big(\hat r_{R',R}^{Q',Q}(\hat x)\big) = 
\,\overline{\!\big(\hat y,\hat r_{V,V}^{R'',R}(\hat r_{R'',R}^{Q'',Q}(\hat s))\,\big)\!} = \,\overline{\!\big(\hat y,\hat r_{V,V}^{Q'',Q}(\hat s)\,\big)\!} =\hat r^{Q',Q}_{V',V}(\hat x)
\eqno £3.7.31.$$

\smallskip
As above, consider a third $V''\in \frak X - \frak Y$ fully normalized in $\F\,,$ 
and  a subgroup $Q''$ of $P$  contained in $N'' = N_P (V'')$ and  strictly containing~$V''\,;$ thus, we have the three $k^*\-$set homomorphisms 
$\hat r^{Q',Q}_{V',V}\,,$ $\hat r^{Q'',Q'}_{V'',V'}$ and $\hat r^{Q'',Q}_{V'',V}$
and we claim that they are compatible with the $k^*\-$compositions, namely 
that we have the following commutative diagram 
$$\matrix{\hat\P^{^\frak Y}\! (Q'',Q')_{V'',V'} \times 
\hat\P^{^\frak Y}\! (Q',Q)_{V',V}&\too &\hat\P^{^\frak Y}\! (Q'',Q)_{V'',V}\cr
\hskip-55pt{\scriptstyle \hat r^{Q'',Q'}_{V'',V'}\times \hat r^{Q',Q}_{V',V}}\hskip10pt\big\downarrow&\phantom{\Big\downarrow}&\big\downarrow\hskip5pt{\scriptstyle \hat r^{Q'',Q}_{V'',V}}\hskip-20pt\cr
\hat\P^{^\frak X}\! (V'',V') \times 
\hat\P^{^\frak X}\! (V',V)&\too &\hat\P^{^\frak X}\! (V'',V)\cr}
\eqno £3.7.32.$$
\smallskip
Indeed, let $\hat x$ and $\hat x'$ be respective elements in 
$\hat\P^{^\frak Y}\! (Q',Q)_{V',V}$ and in $\hat\P^{^\frak Y}\! 
(Q'',Q')_{V'',V'}\,;$ we actually may assume that 
$$\hat x = \hat r^{N',N}_{Q',Q} (\hat y)\.\hat s\qq \hat x' =
 \hat r^{N'',N'}_{Q'',Q'} (\hat y')\.\hat s'
\eqno £3.7.33\phantom{.}$$
where $\hat y$ and $\hat y'$ are suitable elements respectively belonging to 
$\hat\P^{^\frak Y}\! (N',N)_{V',V}$ and $\hat\P^{^\frak Y}\! (N'',N')_{V'',V'}\,,$
and where, denoting by $y$ and $y'$ their images in $\P$ and setting 
$$R = \big(\pi_{N,N'}(y^{-1})\big)(Q')\qq R' = \big(\pi_{N',N''}(y'^{-1})\big)(Q'')
\eqno £3.7.34,$$
$\hat s$ and $\hat s'$ are suitable elements respectively belonging to 
$\hat\P^{^\frak Y}\! (R,Q)_{V,V}$ and to~$\hat\P^{^\frak Y}\!
(R',Q')_{V',V'}\,.$ Then, setting 
$$R'' = \big(\pi_{N,N'}(y^{-1})\big)(R') = \big(\pi_{N,N''}(y'\.y)^{-1}\big)(Q'')
\eqno £3.7.35,$$ 
we clearly have
$$\eqalign{\hat x'\.\hat x &= \big(\hat r^{N'',N'}_{Q'',R'} (\hat y')\.\hat s'\big)\.\big(\hat r^{N',N}_{Q',R} (\hat y)\.\hat s\big)\cr
&= \hat r_{Q'',R''}^{N'',N} (\hat y'\.\hat y)\.\big(\hat r_{R'',R'}^{N,N'}(\hat y^{-1})
\.\hat s'\.\hat r^{N',N}_{Q',R} (\hat y)\big)\.\hat s\cr} 
\eqno £3.7.36.$$
Hence, setting $\hat s'' = \hat r_{R'',R'}^{N,N'}(\hat y^{-1})
\.\hat s'\.\hat r^{N',N}_{Q',R} (\hat y)\,,$ we get (cf.~£3.7.27)
$$\hat r^{Q'',Q}_{V'',V} (\hat x'\.\hat x) = \,\overline{\!\big(\hat y'\.\hat y,\hat r_{V,V}^{R'',Q}(\hat s''\.\hat s)\,\big)\!}
\eqno £3.7.37.$$

\smallskip
On the other hand, from equalities~£3.7.33 we obtain (cf.~£3.7.27)
$$\hat r^{Q',Q}_{V',V} (\hat x) =   
\,\overline{\!\big(\hat y,\hat r_{V,V}^{R,Q}(\hat s)\,\big)\!}\qq 
\hat r^{Q'',Q'}_{V'',V'} (\hat x') =   
\,\overline{\!\big(\hat y',\hat r_{V',V'}^{R',Q'}(\hat s')\,\big)\!}
\eqno £3.7.38;$$
but, according to our definition in~£3.7.19, we get
$$\eqalign{c^{_\frak X}_{V'',V',V} 
&\Big(\,\overline{\!\big(\hat y',\hat r_{V',V'}^{R',Q'}(\hat s')\big)\!}\,,
\,\overline{\!\big(\hat y,\hat r_{V,V}^{R,Q}(\hat s)\big)\!}\,\Big)\cr
& = \,\overline{\!\big(\hat y'\.\hat y,(\hat\frak h_y)^{-1}\big(\hat r_{V',V'}^{R',Q'}(\hat s')\big)\.\hat r_{V,V}^{R,Q}(\hat s)\big)\!}\,\cr}
\eqno £3.7.39\phantom{.}$$
and we claim that we have
$(\hat\frak h_y)^{-1}\big(\hat r_{V',V'}^{R',Q'}(\hat s')\big) 
= \hat r_{V,V}^{R'',R}(\hat s'')$ which will force (cf. £3.7.37)
$$\eqalign{c^{_\frak X}_{V'',V',V} &\Big(\,\overline{\!
\big(\hat y',\hat r_{V',V'}^{R',Q'}(\hat s')\big)\!}\,,\,
\overline{\!\big(\hat y,\hat r_{V,V}^{R,Q}(\hat s)\big)\!}\,\Big)\cr
 &= \,\overline{\!\big(\hat y'\.\hat y,\hat r_{V,V}^{R'',Q}(\hat s''\.\hat s)\,\big)\!}
= \hat r^{Q'',Q}_{V'',V} (\hat x'\.\hat x)\cr}
\eqno £3.7.40$$
completing the proof of the commutativity of diagram £3.7.32.

\smallskip
Denoting by $\varphi'$ the image of $\hat\tau^{_\frak Y}_{N',R'} (1)\.\hat s'$ in 
$\big(N_\F (V')\big) (N',Q')$ (cf.~£3.7.9) and employing the terminology in [5,~5.15],
we argue by induction on the {\it length\/} $\ell (\varphi')$ of~$\varphi'\,;$
if $\ell (\varphi') = 0$ we have $\varphi' = \sigma'\circ \iota_{Q'}^{N'}$ for  
$\sigma'\in \big(N_\F (V')\big) (N')$ [5,~Corollary~5.14] and therefore we get 
$\hat\tau^{_\frak Y}_{N',R'} (1)\.\hat s' = \hat t'\.\hat\tau^{_{\frak Y}}_{N',Q'} (1)$ for a suitable $\hat t'\in \hat\P^{^{\frak Y}}\! (N')_{V'}\,,$ so that we obtain (cf.~£3.7.18)
$$(\hat\frak h_y)^{-1}\big(\hat r_{V',V'}^{R',Q'}(\hat s')\big) =  
\hat r_{V}^{N}\big(\hat g_y (\hat t')\big) = \hat r_{V}^{N}(\hat y^{-1} \.\hat t'\.\hat y) =
\hat r_{V,V}^{R'',R}(\hat s'')
\eqno £3.7.41.$$
Otherwise, we have [5,~5.15.1]
$$\varphi' = \iota_{T'}^{N'}\circ\tau'\circ \eta'\qq 
\ell (\iota_{T'}^{N'}\circ\eta') = \ell (\varphi') -1
\eqno £3.7.42\phantom{.}$$
for  some $T'$ in $N_\frak Y (V')\,,$ some  $\eta'$ in $ \big(N_\F (V')\big) (T',Q')$ and some  
$\tau'$ in\break $\big(N_\F (V')\big) (T')\,,$  
and therefore we get $\hat s' = \hat\tau^{_{\frak Y}}_{N',T'} (1)\.\hat t'\.\hat u'$
for suitable elements $\hat t'\in \hat\P^{^{\frak Y}}\! (T')_{V'}$ and 
$\hat u'\in \hat\P^{^{\frak Y}}\! (T',Q')_{V',V'}$ respectively lifting $\tau'$ and $\eta'\,;$
hence, we obtain 
$$\hat r_{V',V'}^{R',Q'}(\hat s') = \hat r_{V'}^{T'}(\hat t')\.\hat r_{V',V'}^{T',Q'}(\hat u') 
\eqno £3.7.43\phantom{.}$$
and therefore we still obtain 
$$(\hat\frak h_y)^{-1}\big(\hat r_{V',V'}^{R',Q'}(\hat s')\big) = 
(\hat\frak h_y)^{-1}\big(\hat r_{V'}^{T'}(\hat t')\big)\.
(\hat\frak h_y)^{-1}\big(\hat r_{V',V'}^{T',Q'}(\hat u')\big) 
\eqno £3.7.44.$$

\smallskip
Then, by the induction hypothesis,  setting $T = \big(\pi_{N,N'}(y^{-1})\big)(T')$
and $\hat u'' = \hat r_{T,T'}^{N,N'}(\hat y^{-1})
\.\hat u'\.\hat r^{N',N}_{Q',R} (\hat y)\,,$ we have $(\hat\frak h_y)^{-1}
\big(\hat r_{V',V'}^{T',Q'}(\hat u')\big) = \hat r_{V,V}^{T,R}(\hat u'')\,;$
moreover, it is quite clear that in~£3.7.18 replacing $N$  by $T$ and $N'$ by $T'$
we still\break
\eject
\noindent
 get  the commutative  diagrams of $k^*\-$groups
$$\matrix{\hat\P^{^\frak Y}\! (T')_{V'}&\too &\hat\P^{^\frak X}\! (V')\cr
\wr\Vert&&\hskip-15pt{\scriptstyle \hat\frak h_x}\hskip5pt\wr\!\Vert\cr
\hat\P^{^\frak Y}\! (T)_V&\too &\hat\P^{^\frak X}\! (V)\cr}\qq 
\matrix{\hat\P^{^\frak Y} (T')_{V'}&\too &\hat\P^{^\frak Y} (T')\cr
\wr\Vert&&\wr\Vert\cr
\hat\P^{^\frak Y} (T)_V&\too &\hat\P^{^\frak Y} (T)\cr}
\eqno £3.7.45$$
and thus, since $\hat t'$ belongs to $\hat\P^{^\frak Y}\! (T')_{V'}\,,$
setting $\hat t'' = \hat r_{T,T'}^{N,N'}(\hat y^{-1})
\.\hat t'\.\hat r^{N',N}_{T',T} (\hat y)$ we still have $(\hat\frak h_y)^{-1}
\big(\hat r_{V'}^{T'}(\hat t')\big) = \hat r_{V}^{T}(\hat t'')\,.$
Finally, it is easy to check  that $\hat r_{V,V}^{R'',R}(\hat s'') = \hat r_{V}^{T}(\hat t'')\.\hat r_{V,V}^{T,R}(\hat u'')\,,$ which completes the proof of our claim.

\smallskip 
We are ready to define the $k^*\-$set $\hat\P^{^{\frak X}}\! (V',V)$ for any pair of
 subgroups $V$ and $V'$ in $\frak X -\frak Y\,;$ we clearly have 
 $N = N_P (V)\not= V$ and it follows from [5,~Proposition~2.7] that there is an 
 $\F\-$morphism $\nu\, \,\colon N\to P$ such that  $\nu (V)$ is fully normalized in~$\F\,;$ moreover, we choose $\hat n\in \hat\P^{^{\frak Y}}\! 
 \big(\nu (N),N\big)$ lifting the $\F\-$isomorphism $\nu_*$ determined by 
 $\nu\,.$ That is to say, we may assume that
\smallskip
\noindent
£3.7.46\quad {\it There is a pair $(N,\hat n)$ formed by a subgroup $N$ of $P$ which strictly contains and normalizes $V\,,$ and by an element $\hat n$ in 
$\hat\P^{^{\frak Y}}\!\big(\nu (N),N\big)$ lifting~$\nu_*$ for a  
$\F\-$morphism $\nu\,\colon N\to P$ such that $\nu (V)$ is fully normalized
 in~$\F\,.$\/}
\smallskip
\noindent
 We denote by $\hat\frak N(V)$ the set of such pairs and often we write $\hat n$
instead of~$(N,\hat n)\,,$ setting ${}^n \!N = \nu (N)\,,$ ${}^n \!V = \nu (V)\,,$ and $\pi_n = \nu_*$ where $n$ is the image of $\hat n$ 
in~$\P\big(\nu (N),N\big)\,.$

\smallskip
For another pair $(\bar N,\bar n)$ in $\hat\frak N(V)\,,$  denoting by 
$\bar\nu\,\colon \bar  N\to P$  the  $\F\-$mor-phism determined by 
$\hat{\bar n}\,,$ setting $M = \langle N,\bar N\rangle$ and considering a new 
$\F\-$mor-phism $\mu\, \colon M\to P$  such that $\mu (V)$ is fully normalized in~$\F\,,$  we can obtain a third pair $(M,\hat m)$ in $\hat\frak N(V)\,;$ then,  
$\hat r_{{}^m\! N, N}^{{}^m M,M}(\hat m)\.\hat n^{-1}$ 
and $\hat r_{{}^m\!\hat N,\hat N}^{{}^m M,M}(\hat m)\.\hat{\bar n}^{-1}$
respectively belong to $\hat\P^{^{\frak Y}}\! ({}^m\! N,{}^n\! N)$
and to $\hat\P^{^{\frak Y}}\!({}^m\!\bar N,{}^{\bar n}\! \bar N\big)\,;$
in particular, since ${}^n  V\,,$ ${}^{\bar n}  V$ and ${}^mV$ are fully
normalized in~$\F\,,$ the $k^*\-$sets $\hat\P^{^{\frak X}}\!({}^mV,{}^n  V) \,,$ 
$\hat\P^{^{\frak X}}\! ({}^m V,{}^{\bar n}  V) $ and 
$\hat\P^{^{\frak X}}\!  ({}^{\bar n}  V,{}^n  V)$ 
have been already defined above, and we consider the element (cf. £3.7.19)
$$\hat g_{\hat{\bar n},\hat n} =
\hat r_{{}^m V,{}^{\bar n}  V}^{{}^m \! \bar N,{}^{\bar n}\! \bar N}
\big(\hat r_{{}^m\!\hat N,\hat N}^{{}^m M,M}(\hat m)\.
\hat{\bar n}^{-1}\big)^{-1}\.\hat r_{{}^m V,{}^{n} V}^{{}^m\! N,{}^{n}\! N}
\big(\hat r_{{}^m\! N, N}^{{}^m M,M}(\hat m)\.\hat n^{-1}\big) 
\eqno £3.7.47\phantom{.}$$
in $\hat\P^{^{\frak X}}\! ({}^{\bar n}  V,{}^n  V)\,,$ which actually does not depend on the choice of $m\,.$

\smallskip
Indeed, for another pair $(M,\hat m')$ in $\frak N(V)$ we have
$$\eqalign{\hat r_{{}^{m'}\! N, N}^{{}^{m'}\! M,M}(\hat m') 
&= \hat r_{{}^{m'}\! N, {}^m\! N}^{{}^{m'}\! M,{}^m\!M}(\hat m'\.\hat m^{-1})
\.\hat r_{{}^{m}\! N, N}^{{}^{m}\! M,M}(\hat m)\cr
\hat r_{{}^{m'}\!\hat  N,\bar  N}^{{}^{m'}\! M,M}(\hat m') 
&= \hat r_{{}^{m'}\!\hat  N, {}^m\!\hat  N}^{{}^{m'}\! M,{}^m\!M}
(\hat m'\.\hat m^{-1})\.\hat r_{{}^{m}\!\hat  N,\bar  N}^{{}^{m}\! M,M}(\hat m)\cr}
\eqno £3.7.48\phantom{.}$$
and therefore it follows from equality~£3.7.29 that we get
$$\eqalign{\hat r_{{}^{m'} V,{}^{ n}  V}^{{}^{m'} \!  N,{}^{ n}\!  N}
&\big(\hat r_{{}^{m'}\! N, N}^{{}^{m'} M,M}(\hat m')\. \hat n^{-1}\big) \cr
&= \hat r_{{}^{m'} V,{}^{ n}  V}^{{}^{m'} \!  N,{}^{ n}\!  N}
\big(\hat r_{{}^{m'}\!  N, {}^m\!  N}^{{}^{m'}\! M,{}^m\!M}(\hat m'\.\hat m^{-1})
\.\hat r_{{}^{m}\!  N,  N}^{{}^{m}\! M,M}(\hat m)\.\hat  n^{-1}\big) \cr
&= \hat r_{{}^{m'} V,{}^m  V}^{{}^{m'}\! M,{}^m\!M}(\hat m'\.\hat m^{-1})
\. \hat r_{{}^{m'} V,{}^{ n}  V}^{{}^{m'} \!  N,{}^{ n}\!  N}
\big(\hat r_{{}^{m}\!  N,  N}^{{}^{m}\! M,M}(\hat m)\. \hat n^{-1}\big)\cr
\hat r_{{}^{m'} V,{}^{\bar n}  V}^{{}^{m'} \! \bar N,{}^{\bar n}\! \bar N}
&\big(\hat r_{{}^{m'}\!\bar N,\bar N}^{{}^{m'} M,M}(\hat m')\.
\hat{\bar n}^{-1}\big) \cr
&= \hat r_{{}^{m'} V,{}^{\bar n}  V}^{{}^{m'} \! \bar N,{}^{\bar n}\! \bar N}
\big(\hat r_{{}^{m'}\!\bar  N, {}^m\!\bar  N}^{{}^{m'}\! M,{}^m\!M}
(\hat m'\.\hat m^{-1})\.\hat r_{{}^{m}\!\bar  N,\bar  N}^{{}^{m}\! M,M}
(\hat m)\.\hat{\bar n}^{-1}\big) \cr
&= \hat r_{{}^{m'} V,{}^m  V}^{{}^{m'}\! M,{}^m\!M}(\hat m'\.\hat m^{-1})
\. \hat r_{{}^{m'} V,{}^{\bar n}  V}^{{}^{m'} \! \bar N,{}^{\bar n}\! \bar N}
\big(\hat r_{{}^{m}\!\bar  N,\bar  N}^{{}^{m}\! M,M}
(\hat m)\.\hat{\hat n}^{-1}\big)\cr}
\eqno £3.7.49,$$
which proves our claim. Similarly, for any triple of pairs $(N,\hat n)\,,$ 
$(\bar N,\hat{\bar n})$ and $(\skew3\bar {\bar  N},\hat{\bar{\bar n}})$  in 
$\hat\frak N(V)\,,$ considering a pair 
$\big(\langle N, \bar N, \skew3\bar {\bar  N}\rangle,\hat m\big)$ in $\hat\frak N(V)\,,$ it follows from equality~£3.7.29 and from the commutativity of diagram~£3.7.32 that
$$\hat g_{\hat{\bar{\bar n}},\hat n}\.\hat g_{\hat{\bar n},\hat n}  = 
\hat g_{\hat{\bar{\bar n}},\hat n}
\eqno £3.7.50.$$ 
Note that if $V$ is fully normalized in $\F$ then  the pair formed by $N = N_P(V)$
 and by the identity element $\hat\imath_N$ in $\hat\P^{^{\frak Y}}\! (N)$ 
 belongs to~$\hat\frak N(V)\,.$

\smallskip
Then, for any pair of subgroups $V$ and $V'$ in $\frak X -\frak Y\,,$  since for any 
$(N,\hat n)\in hat\frak N(V)$ and any $(N',\hat n')\in \hat\frak N(V')$ the $k^*\-$set 
$\hat\P^{^{\frak X}}\!  ({}^{n'}\! V',{}^{n} V)$ is already defined,
we denote by~$\hat\P^{^{\frak X}}\!  (V',V)$ the $k^*\-$subset of the product
$$ \prod_{\hat n\in\ \hat\frak N(V)}\,\prod_{\hat n'\in \hat\frak N(V')} \hat\P^{^{\frak X}}\!  ({}^{n'}\! V',{}^{n} V)
\eqno £3.7.51\phantom{.}$$
formed by the families $\{\hat x_{\hat n',\hat n}\}_{\hat n\in \hat\frak N(V),\hat n'\in \hat\frak N(V')}$ fulfilling
$$\hat g_{\hat{\bar n}',\hat n'}\.\hat x_{\hat n',\hat n} = \hat x_{\hat{\bar n}',\hat{\bar n}}\. \hat g_{\hat{\bar n},\hat n} 
\eqno £3.7.52.$$ 
In other words, the set $\hat\P^{^{\frak X}}\!  (V',V)$ is the {\it inverse limit\/} of the family formed by the $k^*\-$sets 
$\hat\P^{^{\frak X}}\! \big({}^{n'}\! V', {}^{n} V\big)$ and by the bijections between them induced by the $\hat\P^{^{\frak X}}\!\-$morphisms 
$\hat g_{\hat{\bar n},\hat n}$ and~$\hat g_{\hat{\bar n}',\hat n'}\,.$

\smallskip
Note that, according to equalities~£3.7.50, the {\it projection map\/} onto the factor labeled by the pair $\big((N,\hat n),(N',\hat n')\big)$ induces a $k^*\-$set isomorphism 
$$\frak n_{\hat n',\hat n} : \hat\P^{^{\frak X}}\! (V',V)\cong 
\hat\P^{^{\frak X}}\!\big({}^{n'}\! V',{}^{n} V\big)
\eqno £3.7.53;$$ 
in particular, if $V$ and $V'$ are fully normalized in $\F\,,$ setting $N = N_P (V)$ and  $N' = N_P (V')\,,$  the pairs $(N,\hat\imath_N)$ and $(N',\hat\imath_{N'})$ respectively belong
 to~$\hat\frak N(V)$ and to $\hat\frak N(V')\,,$ and therefore we have a {\it canonical\/} bijection
 $$\frak n_{\hat\imath_{N'},\hat\imath_N} : \hat\P^{^{\frak X}}\! (V',V)\cong 
 \hat\P^{^{\frak X}}\!\big(\,{}^{\hat\imath_{N'}}\! V',{}^{\hat\imath_N} V\big)
\eqno £3.7.54,$$ 
so that our notation is coherent. Moreover, we have an obvious map
$$ \hat\P^{^{\frak X}}\! (V',V)\too \P (V',V)
\eqno £3.7.55\phantom{.}$$
and, for any $u\in\T_{\! P} (V',V)$ and a suitable  pair 
$\big((N,\hat n),(N',\hat n')\big)\,,$ we may assume that $u$ belongs to 
$\T_{\! P} (N',N)$ too; then,  we consider the map 
$$\hat\tau_{V',V}^{_{\frak X}}\! :  \T_{\! P} (V',V)\too
\hat\P^{^{\frak X}}\! (V',V)
\eqno £3.7.56\phantom{.}$$
 determined~by
$$\frak n_{\hat n',\hat n}\big(\hat\tau_{V',V}^{_{\frak X}} (u)\big) = 
\hat r_{{}^{n'}\!V',{}^n V}^{{}^{n'}\!\!N',{}^n\! N}
\big(\hat n'\.\hat\tau_{N',\hat N}^{_{\frak Y}}(u)\.\hat n^{-1}\big)
\eqno £3.7.57,$$
which does not depend on our choice.

\smallskip
Analogously,  for any pair of subgroups  $Q$ and $Q'$ of $P$ respectively normalizing and strictly containing $V$ and $V'\,,$ we can define  an injective 
$k^*\-$set homomorphism
$$\hat r^{Q',Q}_{V',V} : \hat\P^{^\frak Y}\!  (Q',Q)_{V',V}\too 
\hat\P^{^\frak X}\!  (V',V)
\eqno £3.7.58\phantom{.}$$
which lifts the restriction map (cf. £3.2.1)
$$r^{Q',Q}_{V',V} : \P (Q',Q)_{V',V}\too \P (V',V)
\eqno £3.7.59\phantom{.}$$
and coincides with the $k^*\-$set homomorphism~£3.7.23 whenever $V$ and $V'$ are fully normalized in $\F\,;$ indeed, it is clear that we have pairs $(Q,\hat n)$ 
in $\hat\frak N(V)$ and $(Q',\hat n')$ in $\hat\frak N(V')\,,$ and then, for any 
$\hat x\in \hat\P^{^\frak Y}\!  (Q',Q)_{V',V}\,,$ we set
$$\frak n_{\hat n',\hat n}\big(\hat r^{Q',Q}_{V',V} (\hat x)\big) =
\hat r_{{}^{n'}\!V',{}^n V}^{{}^{n'}\!\!Q',{}^n\! Q}
\big(\hat n'\.\hat x\.\hat n^{-1}\big)
\eqno £3.7.60,$$
which does not depend on our choices. Moreover, it is easily checked that 
equality~£3.7.29 still holds in this general situation.

\smallskip
On the other hand, for any $V''\in \frak X -\frak Y\,,$ the {\it $k^*\-$composition map\/} defined in~£3.7.19 --- and just noted $\.$ from now on --- can be extended to a new {\it $k^*\-$composition map\/}
$$\hat\P^{^{\frak X}}\! (V'',V')\times \hat\P^{^{\frak X}}\! (V',V)\too 
\hat\P^{^{\frak X}}\! (V'',V)
\eqno £3.7.61\phantom{.}$$
sending $(\hat x',\hat x)\in \hat\P^{^{\frak X}}\! (V'',V')\times 
\hat\P^{^{\frak X}}\! (V',V)$ to 
$$\hat x'\. \hat x = (\frak n_{\hat n'',\hat n})^{-1}
\big(\frak n_{\hat n'',\hat n'}(\hat x')\.\frak n_{\hat n',\hat n}(\hat x)\big)
\eqno £3.7.62\phantom{.}$$
\eject
\noindent
for a choice of $(N,\hat n)$ in $\hat\frak N(V)\,,$ of $(N',\hat n')$ in $\hat\frak N(V')$ and of $(N'',\hat n'')$ in $\hat\frak N(V'')\,.$ This  {\it $k^*\-$composition map\/} does not depend on our choice; indeed, for another choice of pairs 
$(\bar N,\hat{\bar n})\in \hat\frak N(V)\,,$ $(\bar N',\hat{\bar n}')\in \hat\frak N(V')$ and $(\bar N'',\hat{\bar n}'')\in \hat\frak N(V'')\,,$
we get (cf.~£3.7.52)
$$\eqalign{ \hat g_{\hat n'',\hat{\bar n}''} 
&\.\frak n_{\hat{\bar n}'',\hat{\bar n}'}(\hat x')\.
\frak n_{\hat{\bar n}',\hat{\bar n}} (\hat x) = 
\frak n_{\hat n'',\hat n'}(\hat x')\. \hat g_{\hat n',\hat{\bar n}}\.
\frak n_{\hat{\bar n}',\hat{\bar n}}(\hat x)\cr 
& = \frak n_{\hat n'',\hat n'}(\hat x')\.\frak n_{\hat n',\hat n}(\hat x)\.
\hat g_{\hat n,\hat{\bar n}} = \frak n_{\hat n'',\hat n}(\hat x'\.\hat x)\.
\hat g_{\hat n,\hat{\bar n}}\cr}
\eqno £3.7.63.$$

\smallskip
In particular, for any triple of subgroups $Q\,,$ $Q'$ and $Q''$ of $P$  respectively normalizing and strictly containing $V\,,$ $V'$ and $V''\,,$ choosing pairs
 $(Q,\hat n)$ in $\hat\frak N(V)\,,$  $(Q',\hat n')$ in $\hat\frak N(V')$ and  
 $(Q'',\hat n'')$ in $\hat\frak N(V'')\,.$ the commutativity of the corresponding 
 diagram~£3.7.32 forces  the commutativity of the analogous diagram in the general situation
$$\matrix{\hat\P^{^\frak Y}\! (Q'',Q')_{V'',V'} \times 
\hat\P^{^\frak Y}\! (Q',Q)_{V',V}&\too &\hat\P^{^\frak Y}\! (Q'',Q)_{V'',V}\cr
\hskip-55pt{\scriptstyle \hat r^{Q'',Q'}_{V'',V'}\times \hat r^{Q',Q}_{V',V}}\hskip10pt\big\downarrow&\phantom{\Big\downarrow}&\big\downarrow\hskip5pt{\scriptstyle \hat r^{Q'',Q}_{V'',V}}\hskip-20pt\cr
\hat\P^{^\frak X}\! (V'',V') \times 
\hat\P^{^\frak X}\! (V',V)&\too &\hat\P^{^\frak X}\! (V'',V)\cr}
\eqno £3.7.64.$$
Finally, for any $V'''\in \frak X -\frak Y$ and any $\hat x''\in \hat\P^{^{\frak X}}\! (V''',V'')\,,$ it follows from~£3.7.21 that
$$(\hat x''\.\hat x')\. \hat x = \hat x''\.(\hat x'\.\hat x)
\eqno £3.7.65.$$

\smallskip
We are ready to complete our construction of the announced {\it regular central  
$k^*\-$extension $\hat\P^{^\frak X}\!$ of $\P^{^\frak X}\,,$  endowed with a lifting $\hat\tau^{_{\frak X}}\,\colon \T^{^{\frak X}}_P \to \hat\P^{^{\frak X}}$ of $\tau^{_{\frak X}}$ fulfilling condition~£3.5.1\/};  we are already assuming that 
$\hat\P$ contains $\hat\P^{^\frak Y}\!$ as a {\it full\/} 
$k^*\-$subcategory over~$\frak Y$ and that $\hat\tau$
extends $\hat\tau^{_{\frak Y}}$.
For any subgroups $V$ in $\frak X -\frak Y$ and $Q$ in $\frak Y$ we define
$$\hat\P^{^{\frak X}}\! (V,Q) = \emptyset\qq 
\hat\P^{^{\frak X}}\! (Q,V) = \bigsqcup_{V'} \hat\P^{^{\frak X}}\!(V',V)
\eqno £3.7.66\phantom{.}$$
where $V'$ runs over the set of subgroups $V'\in \frak X -\frak Y$ contained in 
$Q$ and the $k^*\-$subset $\hat\P^{^{\frak X}}\!(V',V)$ of 
$\hat\P^{^{\frak X}}\! (Q,V)$ coincides with the converse image 
of the subset $\tau_{Q,V'}(1)\.\P (V',V)$ in $\P (Q,V)\,;$ moreover, any
$u\in \T_P (Q,V)$ also belongs to $\T_P (uVu^{-1},V)$ and we define 
$\hat\tau_{Q,V}^{_{\frak X}}\! (u)$ as the element 
$\hat\tau_{uVu^{-1},V}^{_{\frak X}}\! (u)$ (cf.~£3.7.56) in the union above.

\smallskip
 In order to define the composition of two $\hat\P^{^{\frak X}}\!\-$morphisms 
$\hat x\,\colon R\to Q$ and $\hat y\,\colon T\to R$ we already may assume that 
$T$ does not belong to $\frak Y\,;$ if $Q$ does not belong to~$\frak Y$ then the composition $\hat x\.\hat y$ is given by the map~£3.7.61;
if $Q$ belongs to $\frak Y$ but $R$ does not  then, setting 
$R' = \varphi (R)$ where $\varphi$ is the image of $\hat x$ in $\F (Q,R)\,,$ 
it follows from definition~£3.7.66 that $\hat x$ is actually an element 
of~$\hat\P^{^{\frak X}}\!(R',R)\,,$ that $\hat y$ is an element of 
$\hat\P^{^{\frak X}}\! (R,T)$ and that the element $\hat x\. \hat y$ defined by the map~£3.7.61 belongs  to~$\hat\P^{^{\frak X}}\!(R',T)\i 
\hat\P^{^{\frak X}}\!(Q,T)\,,$ so that we can define the composition of $\hat x$ and $\hat y$ by this element $\hat x\. \hat y\,.$
Finally, assume that $R$ belongs to $\frak Y$ and, denoting by $\psi$ the image of
$\hat y$ in~$\F (R,T)\,,$ consider the  subgroups $T' = \psi (T)$ of $R$ and 
$T'' = \varphi (T')$ of $Q\,;$ then, it follows again from definition~£3.7.66 
that $\hat y$ is actually an element of $\hat\P^{^{\frak X}}\! (T',T)\,;$
moreover, setting $\bar R = N_R (T')$ and $\bar Q = N_Q (T'')\,,$ it is clear that 
$\hat r_{\bar Q,\bar R}^{Q,R} (\hat x)$ belongs to~$\hat\P^{^{\frak Y}}\! 
(\bar Q,\bar R)$ (cf.~£3.6) and we can define (cf.~£3.7.58 and~£3.7.61)
$$\hat x\. \hat y = \hat r_{T'',T'}^{\bar Q,\bar R}\big(\hat r_{\bar Q,\bar R}^{Q,R} (\hat x)\big)\. \hat y
\eqno £3.7.67.$$

\smallskip
This composition is clearly compatible with the action of $k^*\,.$ Moreover, for a third $\hat\P^{^{\frak X}}\!\-$morphism $\hat z\,\colon V\to T$ we claim that
$$(\hat x\. \hat y)\. \hat z = \hat x \. (\hat y\. \hat z)
\eqno £3.7.68.$$
Once again, we may assume that $V$ does not belong to $\frak Y\,;$  if $Q$ does not belong to~$\frak Y$ then this equality follows from equality~£3.7.65;
if $Q$ belongs to~$\frak Y$ but $R$ does not  then $\hat x$ is actually an element of~$\hat\P^{^{\frak X}}\!(R',R)$ and this equality follows again from equality~£3.7.65. From now on, assume that $R$ belongs to~$\frak Y\,;$
then, if $T\in \frak Y\,,$ denoting by  $\eta$ the image of
$\hat z$ in~$\F (T,V)\,,$ considering the subgroups $V' = \eta (V)$ of $T\,,$
$V'' = \psi (V')$ and $V''' = \varphi (V'')$ and setting  
$\skew3\bar{\bar T} = N_T (V')\,,$ $\skew3\bar{\bar R} = N_R (V'')$ and 
$\skew3\bar{\bar Q} = N_Q (V''')\,,$ then we have (cf.~£3.7.67)
$$(\hat x\. \hat y)\. \hat z = 
\Big(\hat r_{V''',V'}^{\skew3\bar{\bar Q} ,\skew3\bar{\bar T} }
\big(\hat r_{\skew3\bar{\bar Q},\skew3\bar{\bar T} }^{Q,T} 
(\hat x\. \hat y)\big)\Big)\. \hat z
\eqno £3.7.69;$$
but, it follows from~£3.6 and from the commutativity of diagram~£3.7.64 that
$$\hat r_{V''',V'}^{\skew3\bar{\bar Q} ,\skew3\bar{\bar T} }
\big(\hat r_{\skew3\bar{\bar Q},\skew3\bar{\bar T} }^{Q,T} 
(\hat x\. \hat y)\big) = \hat r_{V''',V''}^{\skew3\bar{\bar Q} ,\skew3\bar{\bar R} }
\big(\hat r_{\skew3\bar{\bar Q},\skew3\bar{\bar R} }^{Q,R} 
(\hat x)\big)\.\hat r_{V'',V'}^{\skew3\bar{\bar R} ,\skew3\bar{\bar T} }
\big(\hat r_{\skew3\bar{\bar R},\skew3\bar{\bar T} }^{R,T} (\hat y)\big)
\eqno £3.7.70;$$
consequently, since $\hat y\. \hat z$ is actually an element 
of~$\hat\P^{^{\frak X}}\!(V'',V)\,,$ it follows from equality~£3.7.65 that
$$\eqalign{(\hat x\. \hat y)\. \hat z &= \hat r_{V''',V''}^{\skew3\bar{\bar Q} ,\skew3\bar{\bar R} }\big(\hat r_{\skew3\bar{\bar Q},\skew3\bar{\bar R} }^{Q,R} 
(\hat x)\big)\.\Big(\hat r_{V'',V'}^{\skew3\bar{\bar R} ,\skew3\bar{\bar T} }
\big(\hat r_{\skew3\bar{\bar R},\skew3\bar{\bar T} }^{R,T} (\hat y)\big)
\. \hat z\Big)\cr
&= \hat r_{V''',V''}^{\skew3\bar{\bar Q} ,\skew3\bar{\bar R} }\big(\hat r_{\skew3\bar{\bar Q},\skew3\bar{\bar R} }^{Q,R} 
(\hat x)\big)\.(\hat y\. \hat z) = \hat x\.(\hat y\. \hat z)\cr}
\eqno £3.7.71.$$
Finally, assume that $T$ does not belong to $\frak Y\,;$ then, we actually have
$V' = T\,,$ $V'' = T'$ and $V''' = T''\,,$ and it follows from~£3.7.65 and~£3.7.67
that
$$\eqalign{(\hat x\. \hat y)\. \hat z &=  \Big(\hat r_{T'',T'}^{\bar Q,\bar R}\big(\hat r_{\bar Q,\bar R}^{Q,R} (\hat x)\big)\. \hat y\Big)\. \hat z = 
\hat r_{V''',V''}^{\bar Q,\bar R}\big(\hat r_{\bar Q,\bar R}^{Q,R} (\hat x)\big)\. (\hat y\. \hat z)\cr
& = \hat x\.(\hat y\. \hat z)\cr}
\eqno £3.7.72.$$

\smallskip
It remains to prove the functoriality of $\hat\tau^{_\frak X}\,;$ that is to say,
for any pair of $\T^{^{\frak X}}_P\-$morphisms $u\,\colon R\to Q$ and 
$v\,\colon T\to R$ we claim that 
$$\hat\tau^{_\frak X}_{Q,T} (uv) = \hat\tau^{_\frak X}_{Q,R} (u)\. 
\hat\tau^{_\frak X}_{R,T} (v)
\eqno £3.7.73;$$
once again, we may assume that $T$ does not belong to $\frak Y\,;$ setting
$T' = vTv^{-1}$ and $T'' = uT'u^{-1}\,,$ it follows easily from our definition and from~£3.7.57 that we have 
$$\eqalign{\hat\tau^{_\frak X}_{Q,T} (uv) &= \hat\tau^{_\frak X}_{T'',T} (uv)
= \hat\tau^{_\frak X}_{T'',T'} (u)\.\hat\tau^{_\frak X}_{T',T} (v)\cr 
\hat\tau^{_\frak X}_{T',T} (v) &= \hat\tau^{_\frak X}_{R,T} (v)\cr}
\eqno £3.7.74;$$
if $R$  does not belong to $\frak Y$ then we have $R = T'$ and, according to our definition, we still have $\hat\tau^{_\frak X}_{T'',T'} (u) = 
\hat\tau^{_\frak X}_{Q,R} (u)\,;$ otherwise,  setting $\bar R = N_R (T')$ and $\bar Q = N_Q (T'')\,,$ it follows from~£3.7.67 and~£3.7.57 that
$$\eqalign{\hat\tau^{_\frak X}_{Q,R} (u)\. \hat\tau^{_\frak X}_{R,T} (v) 
&= \hat r_{T'',T'}^{\bar Q,\bar R}\Big(\hat r_{\bar Q,\bar R}^{Q,R} 
\big(\hat\tau^{_\frak Y}_{Q,R} (u)\big)\Big)\. \hat\tau^{_\frak X}_{R,T} (v)\cr
&= \hat r_{T'',T'}^{\bar Q,\bar R}
\big(\hat\tau^{_\frak Y}_{\bar Q,\bar R} (u)\big)\. \hat\tau^{_\frak X}_{T',T} (v)
= \hat\tau^{_\frak X}_{T'',T'} (u)\.\hat\tau^{_\frak X}_{T',T} (v) \cr}
\eqno £3.7.75.$$

\smallskip
In order to prove the uniqueness of $\hat\P^{^{\frak X}}\,,$ let $\widehat{\P^{^{\frak X}}}$ be
another {\it regular central $k^*\-$extension\/} of $\P^{^\frak X}\,,$ endowed with a functor
$\widehat{\tau^{^\frak X}}  \colon \T^{^\frak X}_P \to \widehat{\P^{^\frak X}}$ fulfilling condition £3.5.1, inducing the  folded Frobenius $P\-$category $(\F,\widehat\aut_{\F^{^{\frak X}}})$ or, equivalently, $(\P,\widehat\aut_{\P^{^{\frak X}}})$. We may assume that $\frak X \not= \{P\}$ and then, choosing a minimal element $U$ in $\frak X$ {\it fully normalized\/} in 
$\F$ and seting
$$\frak Y = \frak X - \{\theta(U)\mid \theta\in \F(P,U)\}
\eqno £3.7.76,$$
we may also assume that $N_\F (U)\not= \F\,.$

\smallskip
In particular, for any group $Q$ in $\frak X$,  denoting by $\frak q_Q \colon \Delta_0\to \P^{^{\frak X}}$ the functor sending $0$ to $Q$, we have
$$\widehat{\P^{^\frak X}} (Q) = \widehat\aut_{\P^{^{\frak X}}} (\frak q_Q) = 
\hat\P^{^{\frak X}}(Q)
\eqno £3.7.77;$$
similarly, for any group $V$ in $\frak X - \frak Y$ fully normalized in $\F$, setting $N = N_P (V)$ and denoting by $\frak n_V \colon \Delta_1 \to P^{^{\frak X}}$ the functor sending $0$ to $V$, $1$ to $N$ and $0\bullet 1$ to $\hat\tau^{^\frak X}_{N,V} (1)$, and by 
$\widehat{\P^{^\frak X}} (N)_V$ and $\hat\P^{^\frak X} (N)_V$ the corresponding stabilizers of $V$ in $\widehat{\P^{^\frak X}} (N)$ and $\hat\P^{^\frak X} (N)$, we have
$$\widehat{\P^{^\frak X}} (N)_V = \widehat\aut_{\P^{^{\frak X}}} (\frak n_V) 
= \hat\P^{^\frak X} (N)_V
\eqno £3.7.78;$$
moreover, $\widehat\aut_{\P^{^{\frak X}}}$ sends the obvious 
$\frak c\frak h^* (\P^{^\frak X})\-$morphism $(\frak n_V,\Delta_1) \to 
(\frak q_V,\Delta_0)$ to the injective restriction from $\widehat{\P^{^\frak X}} (N)_V = \hat\P^{^\frak X} (N)_V$ to $\widehat{\P^{^\frak X}} (V)  = 
\hat\P^{^{\frak X}}(V)$.

\smallskip
Arguing by induction on $\vert \frak X \vert$ we may assume that we have an equivalence of categories $\frak f^{^\frak Y} \colon \widehat{\P^{^\frak Y}} \to \hat\P^{^\frak Y} $ inducing the identity on $\widehat{\P^{^\frak Y}} (Q)  = \hat\P^{^{\frak Y}}(Q)$ for any group 
$Q$ in $\frak Y$ and fulfilling $\frak f^{^\frak Y} \circ \widehat{\tau^{^\frak Y}}   = \hat\tau^{^\frak Y}  $. We will extend $\frak f^{^\frak Y}$ to a functor 
$\frak f^{\frak X} \colon \widehat{\P^{^\frak X}} \to \hat\P^{^\frak X}$  inducing the identity on $\widehat{\P^{^\frak Y}} (Q)  = \hat\P^{^{\frak Y}}(Q)$ for any group 
$Q$ in~$\frak X$ and fulfilling $\frak f^{^\frak X} \circ \widehat{\tau^{^\frak X}}   = \hat\tau^{^\frak X}\,;$ for any pair of groups $V$ and $V'$ in $\frak X - \frak Y$ fully normalized in $\F$, any $\hat y\in \widehat{\P^{^\frak Y}} (N',N)_{V',V}$ where $N' = N_P (V')$ and 
$N = N_P (V)\,,$ and any $\hat s\in \widehat{\P^{^\frak X}} (V)$, we define
$$\frak f^{^\frak X}  \big(\widehat{r{^{_\frak X}}}^{N',N}_{V',V} (\hat y).\hat s\big) =
\hat r{^{_\frak X}}^{N',N}_{V',V} \big( \frak f^{^\frak Y}(\hat y)\big). \hat s
\eqno £3.7.79\,;$$
the definition is correct since for any $\hat t\in \hat\P^{^\frak Y} (N)_V$ we have
$$\eqalign{\frak f^{^\frak X}  \big(\widehat{r{^{_\frak X}}}^{N',N}_{V',V} (\hat y.\hat t).
&(\widehat{r{^{_\frak X}}}^{N,N}_{V,V} (\hat t^{-1})\.\hat s)\big) = \hat r{^{^\frak X}}^{N',N}_{V',V} \big( \frak f^{^\frak Y}(\hat y.\hat t)\big).(\widehat{r{^{_\frak X}}}^{N,N}_{V,V} (\hat t^{-1}). \hat s)\cr
&= \hat r{^{^\frak X}}^{N',N}_{V',V} \big( \frak f^{^\frak Y}(\hat y). \hat t\big).
.(\widehat{r{^{_\frak X}}}^{N,N}_{V,V} (\hat t^{-1}). \hat s)\cr
&=  \hat r{^{^\frak X}}^{N',N}_{V',V} \big( \frak f^{^\frak Y}(\hat y)\big).\hat s\cr}
\eqno £3.7.80.$$

\smallskip
It follows from Lemma £3.3 that $\frak f^{^\frak X}$ induces a bijection from 
$\widehat{\P^{^\frak X}} (V',V)$ onto $\hat\P^{^\frak X} (V',V)\,;$ moreover, if $V''$ is a third group in $\frak X - \frak Y$ fully normalized in $\F\,,$ setting $N'' = N_P (V'')$ and considering $\widehat y'\in \widehat{\P^{^\frak Y}} (N'',N')_{V'',V'}$ and 
$\widehat s'\in \widehat{\P^{^\frak X}}(V')$, it follows from [5, Condition 2.8.2] that 
$\widehat s' = \widehat{r{^{_\frak X}}}^{N',N'}_{V',V'} (\widehat z')$ for some 
$\widehat z'\in \widehat{\P^{^\frak Y} }(N')$ and therefore we get
$$\eqalign{\frak f^{^\frak X} \big(\widehat{r{^{_\frak X}}}^{N'',N'}_{V'',V'} (\widehat y').&\widehat s'  .\widehat{r{^{_\frak X}}}^{N',N}_{V',V} (\widehat y).\widehat s\big) 
= \frak f^{^\frak X} \big(\widehat{r{^{_\frak X}}}^{N'',N}_{V'',V} (\widehat y'.\widehat z' . \widehat y).\widehat s\big)\cr
&={ \hat r{^{_\frak X}}}^{N'',N}_{V'',V}\big(\frak f^{^\frak Y}(\widehat y'.\widehat z' . \widehat y)\big ) .\hat s\cr
&= \hat r{^{_\frak X}}^{N'',N'}_{V'',V'}\big(\frak f^{^\frak Y}(\widehat y')\big) . \hat s'. 
\hat r{^{_\frak X}}^{N',N}_{V',V}\big(\frak f^{^\frak Y}(\widehat y)\big) . \hat s\cr
&= \frak f^{^\frak X} \big(\widehat{r{^{_\frak X}}}^{N'',N'}_{V'',V'} (\widehat y').\widehat s'\big) . 
\frak f^{^\frak X} \big( \widehat{r{^{_\frak X}}}^{N',N}_{V',V} (\widehat y).\widehat s\big)\cr}
\eqno £3.7.81.$$:

\smallskip
In particular, for any group $V$ in $\frak X - \frak Y\,,$ as in £3.7.46 we can define an analogous set $\widehat{\frak N} (V)$ of pairs  $(N,\widehat n)$ formed by a subgroup $N$ of $P$ which strictly contains and normalizes $V\,,$ and by an 
$\widehat{\P^{^\frak Y} }\-$isomorphism $\widehat n$ from $N$ such that ${}^n V\,,$ where $n$ is the image of $\widehat n$ in $N\,,$
is fully normalized in $\F\,;$ similarly, for any pair of elements  $(N,\widehat n)$ and $(\bar N,\widehat{\bar n})$ in $\widehat{\frak N}(V)$, we can define an element 
$\widehat g_{\widehat{\bar n},\widehat n}$ in 
$\widehat{\P^{^\frak X}} ({}^{\bar n}V,{}^nV)$ analogous to the element $\hat g_{\hat{\bar n},n}\in \hat\P^{^\frak X} ({}^{\bar n}V,{}^nV)$ defined in £3.7.47 above and clearly get 
$\frak f^{^\frak X} (\widehat g_{\widehat{\bar n},\widehat n}) = \hat g_{\hat{\bar n},n}\,.$ Then, for any group $V'$ in $\frak X - \frak Y\,,$ we have an obvious bijection from 
$\widehat{\P^{^\frak X}} (V',V)$ onto the $k^*\-$subset of the product
$$ \prod_{\widehat n\in \widehat{\frak N}(V)}\,\prod_{\widehat n'\in \widehat{\frak N}(V')} \widehat{\P^{^{\frak X}}}\!  ({}^{n'}\! V',{}^{n} V)
\eqno £3.7.82\phantom{.}$$
formed by the families $\{\widehat x_{\widehat n',\widehat n}\}_{\widehat n\in \widehat{\frak N}(V),\widehat n'\in \widehat{\frak N}(V')}$ fulfilling
$$\widehat g_{\widehat{\bar n}',\widehat n'}\.\widehat x_{\widehat n',\widehat n} = 
\widehat x_{\widehat{\bar n}',\widehat{\bar n}}\. \widehat g_{\widehat{\bar n},\widehat n} 
\eqno £3.7.83;$$ 
 hence, $\frak f^{^\frak X}$ can be extended to a bijection from 
 $\widehat{\P^{^\frak X}} (V',V)$ onto $\hat{\P^{^\frak X}} (V',V)\,.$ At present, it is quite clear that $\frak f^{\frak X}$ can be extended to an equivalence of categories from 
 $\widehat{\P^{\frak X}}$ onto $\hat\P^{\frak X}\,.$ We are done.
 
 \bigskip
\noindent
{\bf Corollary~£3.8.} {\it Let $G$ be a finite group, $b$ a block of $G$ and $P$
a defect group of $b\,.$ There is  a {\it regular central $k^*\-$extension\/} 
 $\hat\F_{\!(b,G)}^{^{\rm sc}}$ of $\F_{\!(b,G)}^{^{\rm sc}}$ admitting
 a {\it $k^*\-$group isomorphism\/}
$$\hat\F_{\!(b,G)}^{^{\rm sc}} (Q)  \cong
\hat N_G(Q,f)/C_G (Q)
\eqno £3.8.1$$
for any $\F_{\!(b,G)}\-$selfcentralizing subgroup $Q$ of $P\,.$\/}

\medskip
\noindent
{\bf Proof:} It is an easy consequence of [5, Theorem~11.32] and Theorem~£3.7.

\bigskip
\noindent
{\bf References}
{\cds
\medskip
\noindent
[1] Andrew Chermak. {\cdt Fusion systems and localities\/}, 
Acta Mathematica, 211(2013), 47-139.
\smallskip
\noindent
[2] Bob Oliver. {\cdt Existence and Uniqueness of Linking Systems: Chermak's proof via obstruction theory\/}, 
Acta Mathematica, 211(2013), 141-175.
\smallskip
\noindent
[3] Sejong Park, {\cdt The gluing problem for some block fusion systems}, Journal of Algebra, 323 (2010), 1690-1697
\smallskip
\noindent
[4] Llu\'\i s Puig, {\cdt Pointed groups and  construction of modules}, Journal of Algebra, 116(1988), 7-129
\smallskip
\noindent
[5] Llu\'\i s Puig, {\cdt ``Frobenius Categories versus Brauer Blocks''\/}, Progress in Math. 274, 2009, Birkh\"auser, Basel
\smallskip
\noindent
[6] Llu\'\i s Puig, {\cdt Ordinary Grothendieck groups of a Frobenius P-category\/}, Algebra Colloquium 18(2011), 1-76
\smallskip
\noindent
[7] Llu\'\i s Puig, {\cdt Existence, uniqueness and functoriality of the perfect
locality over a Frobenius P-category\/}, Algebra Colloquium, 23(2016), 541-622
\smallskip
\noindent
[8] Llu\'\i s Puig, {\cdt A correction to the uniqueness of a partial perfect
locality over a Frobenius P-category\/}, Algebra Colloquium,  26(2019), 541-559}

\end